\input{amssymb.sty}
\documentclass[11pt]{amsart}
\usepackage[all]{xy}
\newdimen\theight
\def\TeXref#1{%
              \leavevmode\vadjust{\setbox0=\hbox{{\tt
                      \quad\quad  {\small \textrm #1}}}%
              \theight=\ht0
              \advance\theight by \lineskip
              \kern -\theight \vbox to
              \theight{\rightline{\rlap{\box0}}%
              \vss}%
              }}%

\theoremstyle{plain}

\newtheorem{thm}{Theorem}[section]
\newtheorem{mainthm}{Theorem}
\newtheorem{maincor}[mainthm]{Corollary}
\newtheorem*{thm*}{Theorem}
\newtheorem{cor}[thm]{Corollary}

\newtheorem*{cor*}{Corollary}

\newtheorem*{conj*}{Conjecture}
\newtheorem*{lemma*}{Lemma}

\newtheorem{lemma}[thm]{Lemma}
\newtheorem*{prop*}{Proposition}

\newtheorem{prop}[thm]{Proposition}

\theoremstyle{definition}

\newtheorem*{defn*}{Definition}
\newtheorem*{rems*}{Remarks}
\newtheorem*{proof*}{Proof}
\newtheorem{prel*}{Preliminaries}
\newtheorem{examples*}{Examples}

\newcommand{\C}{\mathbb{C}}
\newcommand{\E}{\widetilde{\mathcal E}}
\newcommand{\F}{\widetilde{\mathcal F}}
\newcommand{\M}{\widetilde{M}}

\newcommand{\cl}{{\mathcal L}}

\def\ga{\gamma}

\def\eps{\varepsilon}
\def\ka{\kappa}

\def\pa{\partial}
\def\spec{\hbox{spec}}
\def\Dom{\operatorname{Dom}}

\def\bx{\bar{x}}

\def\Ker{\operatorname{Ker}\,}
\def\Im{\operatorname{Im}\,}

\def\sup{\operatorname{sup}}
\def\inf{\operatorname{inf}}

\def\C{\mathbb C}
\def\R{\mathbb R}

\def\E{\mathcal E}
\def\K{\mathcal K}

\def\N{\mathbb N}
\def\B{\mathcal B}

\def\Tr{\operatorname{Tr}}

\def\A{{\mathcal A}} \def\B{{\mathcal B}}
 
\def\F{{\mathcal F}}

\def\tr{\operatorname{tr}}

\def\id{\operatorname{id}}
\def\H{\mathcal H}

\def\E{\mathcal E}
\def\U{{\mathcal U}}
\def\<{\langle}
\def\>{\rangle}

\newcommand{\supp}{\operatorname{supp}}

\newcommand{\nc}{\newcommand}
\nc{\nt}{\newtheorem}
\nc{\gf}[2]{\genfrac{}{}{0pt}{}{#1}{#2}}
\nc{\mb}[1]{{\mbox{$ #1 $}}}
\nc{\real}{{\mathbb R}}
\nc{\comp}{{\mathbb C}}
\nc{\ints}{{\mathbb Z}}
\nc{\Ltoo}{\mb{L^2({\mathbf H})}}
\nc{\rtoo}{\mb{{\mathbf R}^2}}
\nc{\slr}{{\mathbf {SL}}(2,\real)}
\nc{\slz}{{\mathbf {SL}}(2,\ints)}
\nc{\su}{{\mathbf {SU}}(1,1)}
\nc{\so}{{\mathbf {SO}}}
\nc{\hyp}{{\mathbb H}}
\nc{\disc}{{\mathbf D}}
\nc{\torus}{{\mathbb T}}

\newcommand{\tE}{\widetilde{\E}}

\newcommand{\cH}{{\mathcal H}}
\nc{\ca}{{\mathcal A}}
\nc{\cag}{{{\mathcal A}^\Gamma}}
\nc{\cg}{{\mathcal G}}
\nc{\chh}{{\mathcal H}}
\nc{\ck}{{\mathcal B}}
\nc{\cm}{{\mathcal M}}
\nc{\cs}{{\mathcal S}}
\nc{\cz}{{\mathcal Z}}
\nc{\sind}{\sigma{\rm -ind}}

\begin{document}

\title [Equivalence of spectral projections  in
semiclassical limit]
{Equivalence of spectral projections  in
semiclassical limit and a vanishing theorem for higher traces in
$K$-theory}
\author{Y. Kordyukov}
\address{Institute of Mathematics, Russian Academy of Sciences, Ufa,
Russia}
\email{yuri@imath.rb.ru}
\author{V. Mathai}
\address{Department of Mathematics, University of Adelaide, Adelaide 5005, Australia}
\email{vmathai@maths.adelaide.edu.au}
\author{M. Shubin}
\address{Department of Mathematics, Northeastern University, Boston,
Mass., USA}
\email{shubin@neu.edu}

\thanks{Yu.K. acknowledges hospitality and support of University of 
Adelaide and Northeastern University as well as partial support from 
the Russian Foundation for Basic Research, Grant 04-01-00190, 
the Australian Research Council and NSF grant DMS-9706038. V.M.
acknowledges support from the Clay Mathematical Institute and the
Australian Research Council. M.S. acknowledges partial support
from NSF grants DMS-9706038 and DMS-0107796.}


\begin{abstract} In this paper, we study a refined $L^2$ version of
the semiclassical approximation  of projectively invariant
elliptic operators with invariant Morse type potentials on
covering spaces of compact manifolds. We work on the level of
spectral projections (and not just their traces) and obtain an
information about classes of these projections in $K$-theory in
the semiclassical limit as the coupling constant $\mu$ goes to
zero. An important corollary is a vanishing theorem for the higher
traces in cyclic cohomology for the spectral projections. This
result is then applied to the quantum Hall effect. We also give a
new proof that there are arbitrarily many gaps in the spectrum
of the operators under consideration in the semiclassical limit.
\end{abstract}
\maketitle

\section*{Introduction}

Let $M$ be a closed Riemannian manifold and $\widetilde M$ be its
universal cover. Let $\omega$  be a closed 2-form on $M$ and $\bf
B$ be its lift to $\M$, so that $\bf B$ is a $\Gamma$-invariant
closed 2-form on $\M$ where $\Gamma$ denotes the fundamental group
of $M$ acting on $\M$ by the deck transformations. We assume that
$\bf B$ is exact. Choose a 1-form $\bf A$ on $\M$ such that $d{\bf
A} = \bf B$. As in geometric quantization we may regard $\bf A$ as
defining a Hermitian connection $\nabla_{\bf A} = d+i{\bf A}$ on
the trivial line bundle $\mathcal L$ over $\M$, whose curvature is
$i\bf B$. Physically we can think of $\bf A$ as the
electromagnetic vector potential for a magnetic field $\bf B$.
Suppose that $\E$ is a Hermitian vector bundle on $M$ and $\tE$
the lift of $\E$ to the universal cover $\widetilde M$. Let
$\widetilde \nabla^\E$ denote a $\Gamma$-invariant Hermitian
connection on $\tE$. Then consider the Hermitian connection
$\nabla = \widetilde\nabla^\E\otimes \id + \id\otimes \nabla_{\bf
A}$ on $\tE\otimes {\mathcal L} = \tE$. It is no longer
$\Gamma$-invariant, but it is invariant under a projective action
of $\Gamma$, as will be explained below.

Using the Riemannian
metric on $\widetilde M$ and the Hermitian metric on $\tE$, consider
the elliptic self-adjoint differential operator
given by
\begin{equation}\label{ham}
H(\mu) = \mu \nabla^*\nabla + B  + \mu^{-1} V,
\end{equation}
where $B$, $V$ are $\Gamma$-invariant self-adjoint endomorphisms of
the  bundle $\tE$, $\mu$ is the coupling constant
and where $V$ satisfies
in addition the following {\em Morse type condition:}

For all $x\in \M$, $V(x) \ge 0$. Moreover, if the matrix $V(x_0)$
is degenerate for some $x_0$ in $\M$, then $V(x_0) = 0$ and there
is a positive constant $c$ such that $V(x)\ge c|x-x_0|^2 I$ for
all $x$ in a neighborhood of $x_0$, where $I$ denotes the identity
endomorphism of $\tE$.

We will also assume that $V$ has at least one zero point. We
remark that all functions $V=|df|^2$, where $|df|$ denotes the
pointwise norm of the differential of a $\Gamma$-invariant Morse
function $f$ on $\M$, are examples of Morse type potentials.

We will analyse the qualitative aspects of
the spectrum of $H(\mu)$ acting on the Hilbert space
\begin{equation}\label{hil1}
\mathfrak{H}=L^2(\M, \tE).
\end{equation}
An important
feature of the elliptic operator $H(\mu)$ is that it commutes with  a
projective $(\Gamma, \sigma)$-action of the fundamental group $\Gamma$.
Here
$\sigma$ denotes the multiplier or $U(1)$-valued 2-cocycle on $\Gamma$
defining this projective action.

Associated to  $H(\mu)$, there is a {\em  model operator} $K(\mu)$
(cf. section 2) which is obtained as a direct sum of quadratic
parts of $H(\mu)$ near the degenerate points of $V$ in a
fundamental domain. It acts on the Hilbert space
\begin{equation}\label{hil2}
\mathfrak{H}_K= L^2({\mathbb R}^n, \C^k)^N,
\end{equation}
where $n$ is the dimension of $M$, $\C^k$ is the fibre of $\E$ and
$N$ denotes the number of zeroes of $V$ that lie in a fundamental
domain.  The model operator $K(\mu)$ is a Hilbert direct sum of
harmonic oscillators. If we take a direct sum of all these
harmonic oscillators over $\M$ (and not only in a fundamental
domain) then we will get another version of the model operator
which has the same spectrum and represents the Hamiltonian for the
crystal obtained with perfectly isolated atoms.

Note that $K(\mu)$ is obtained by a simple scaling from the
operator $K=K(1)$ and has a discrete spectrum independent of
$\mu$. Therefore the spacing of its eigenvalues is bounded below.
Then we have the following

\begin{mainthm}[Existence of spectral gaps]
\label{main0} In the notation above, let $V$ be a Morse type
endomorphism of the vector bundle $\tE$ over $\M$ and $H(\mu)$ as
in (\ref{ham}) be the elliptic self-adjoint operator acting on the
Hilbert space $\mathfrak H$ as in  (\ref{hil1}). If $[a, b]$ is an
interval in $\R$ that does not intersect the spectrum of the model
operator $K$ acting on the Hilbert space ${\mathfrak H}_K$, cf.
(\ref{hil2}), then there exists $\mu_0 >0$ such that for all $\mu
\in (0, \mu_0)$, the interval $[a, b]$ also does not intersect the
spectrum of $H(\mu)$.
There exists  arbitrarily
large number of gaps in the spectrum of
$H(\mu)$ provided the coupling constant $\mu$ is sufficiently small.
\end{mainthm}

The physical explanation for the appearance of gaps in the spectrum
$H(\mu)$ is that the potential wells get deeper as $\mu\to 0$
and the atoms get (asymptotically) isolated, so that the energy
levels of $H(\mu)$ are approximated by those of the
corresponding model operator $K$.

The special case of Theorem~\ref{main0} in the absence of a
magnetic field was established in \cite{Sh} using a different
method, and the special case of this result in the presence of a
magnetic field but in the scalar case was established in \cite{MS}
using the same method as in \cite{Sh}.

Theorem~\ref{main1} below is a significant refinement of
Theorem~\ref{main0} above and  will be established using a
refinement of the $L^2$ semiclassical asymptotics used there. We
first set some notation.

If $\cH$ is a Hilbert space, then ${\mathcal K}(\cH)$ denotes the
algebra of compact operators in $\cH$, and $\mathcal{K}=
\mathcal{K}(\ell^2(\N))$, where $\N=\{1,2,3,\ldots\}$.

Let $\mathcal A$ be a unital $*$-algebra with the unit $1_A$, and let
${\rm Proj}(\mathcal A)$ be its set of self-adjoint projections.
Two
projections $P, Q \in {\rm Proj}(\mathcal A)$ are said to be {\it
Murray-von Neumann equivalent} if there is an element $V\in
\mathcal A$ such that $P= V^*V$ and $Q= VV^*$.
Denote $M_n(\C)\otimes \mathcal A = M_n( \mathcal A )$,
where $M_n(\C)$ denotes the square matrices of size
$n$ over $\C$.
Then $M_n( \mathcal A )$ is also a
$*$-algebra. Let $M_\infty( \mathcal A ) = \lim_{n\to\infty}
M_n( \mathcal A )$ be the direct limit of the embeddings of $M_n(
\mathcal A )$ in $M_{n+1}( \mathcal A )$ given by $A\to
\left(\begin{array}{cc}A & 0\\0 & 0\end{array}\right)$. Let
$V(\mathcal A) = {\rm Proj}(M_\infty( \mathcal A ))/\sim\;$ denote
the Murray-von Neumann equivalence classes of projections in
$M_\infty(\mathcal A )$.
Then $V(\mathcal A) $ is
an Abelian  semi-group under with the operation induced by the direct sum,
and the associated Abelian group is called the Grothendieck group 
$K_0(\mathcal A)$.

The homomorphism $\pi:\mathbb C \to \mathcal A$ given by $\lambda
\mapsto \lambda\cdot 1_A$ induces a homomorphism
$\pi_*:K_0(\mathbb C )\cong {\mathbb Z} \to K_0(\mathcal A)$. Then
the reduced $K$-group $\tilde{K}_0(\mathcal A) $ is defined as
$\tilde{K}_0(\mathcal A)= \operatorname{coker} \pi_*\cong
K_0(\mathcal A)/{\mathbb Z}$.

Suppose that  $\mathcal A$ is a non-unital $*$-algebra. Let
$\widetilde{\mathcal A} = \left\{(a, \lambda): a\in \mathcal A,
\;\lambda\in \mathbb C\right\}$. Then $\widetilde{\mathcal A}$ is
a unital $*$-algebra containing $\mathcal A$, with product given
by $(a, \lambda)(b, \mu) = (ab + \lambda b + \mu a, \lambda\mu)$.
By definition, the $K$-group $K_0(\mathcal A)$ is the reduced
$K$-group $\tilde{K}_0(\widetilde{\mathcal A})$ of
$\widetilde{\mathcal A}$.

Recall that the Morita invariance of $K$-theory asserts that there
is a natural isomorphism $K_0(\mathcal A ) \cong K_0(M_n( \mathcal
A ) ) $ which is induced by the standard algebra homomorphism
$\mathcal A \to M_n( \mathcal A )$ which maps $a\in \mathcal A$ to
a matrix with the left-upper corner matrix element $a$, the rest
matrix elements being $0$. If $\mathcal A$ is a $C^*$-algebra,
then the Morita invariance of $K$-theory asserts that there is a
natural isomorphism $K_0(\mathcal A ) \cong K_0(\mathcal A \otimes
\mathcal K) $ which is induced by the similar algebra homomorphism
map $\mathcal A \to \mathcal A \otimes \mathcal K.$

Finally, we denote by $C^*_r(\Gamma,\bar\sigma)$ the reduced
twisted group $C^*$-algebra of the group $\Gamma$. We will assume
that the algebra $C^*_r(\Gamma,\bar\sigma)$ acts on
$\ell^2(\Gamma)$ by left twisted convolutions.

\begin{mainthm}[Semiclassical vanishing theorem in $K$-theory for spectral
projections]\label{main1} In the notation above, let $V$ be a
Morse type endomorphism of the vector bundle $\tE$ over $\M$ and
$H(\mu)$ as in (\ref{ham}) be the elliptic self-adjoint operator
acting on the Hilbert space $\mathfrak H$ as in  (\ref{hil1}). Let
$\lambda\in \C$ be such that $\lambda$ is not in the spectrum of
the model operator $K$ acting on the Hilbert space ${\mathfrak
H}_K$, cf. (\ref{hil2}). Let $E(\lambda) = {\chi}_{(-\infty,
\lambda]}(H(\mu))$  and  $E^0(\lambda) = \chi_{(-\infty,
\lambda]}(K(\mu))$ denote the spectral projections.
\begin{enumerate}

\item There exists a $(\Gamma, \sigma)$-equivariant
isometry $U:{\mathfrak H}\to \ell^2(\Gamma)\otimes {\mathfrak
H}_K$ (see section 2) and a constant $\mu_0 >0$ such that for all
$\mu\in (0, \mu_0)$, the spectral projections $U E(\lambda) U^*$
and $\id\otimes E^0(\lambda)$ are in $C^*_r(\Gamma,
\bar\sigma)\otimes \mathcal K({\mathfrak H}_K)$ and are Murray-von
Neumann equivalent in $C^*_r(\Gamma, \bar\sigma)\otimes \mathcal
K({\mathfrak H}_K)$. In particular,
\begin{align}
[U E(\lambda) U^*]&=[\id\otimes E^0(\lambda)] \in
K_0(C^*_r(\Gamma, \bar\sigma)\otimes \K({\mathfrak H}_K))\cong
K_0(C^*_r(\Gamma, \bar\sigma));\label{e:1} \\ [ E(\lambda) ] &= 0
\in \widetilde K_0(C^*_r(\Gamma, \bar\sigma)). \label{e:11}
\end{align}

\item There is a smooth subalgebra $\B(\Gamma,
\sigma)$  of $C^*_r(\Gamma, \bar\sigma)\otimes \K({\mathfrak
H}_K)$, cf. section 3,
such that the spectral projections $U E(\lambda) U^*$ and $\id
\otimes E^0(\lambda)$ are in $\B(\Gamma, \sigma)$ and are also
Murray-von Neumann equivalent in $\B(\Gamma, \sigma)$.  That is,
for all $\mu\in (0, \mu_0)$, one has $$[U E(\lambda)
U^*]=[\id\otimes E^0(\lambda)] \in K_0(\B(\Gamma, \sigma)).$$
\end{enumerate}
\end{mainthm}

Let $\Tr_{\Gamma}$ denote the trace on $C^*_r(\Gamma,
\bar\sigma)\otimes \K({\mathfrak H}_K)$, which is the tensor
product of the canonical finite trace $\tr_{\Gamma}$ on
$C^*_r(\Gamma, \bar\sigma)$ and the standard trace
$\Tr$ on $\K({\mathfrak H}_K)$. As an immediate consequence of
Theorem~\ref{main1}, we get the following

\begin{maincor}[Semiclassical asymptotics of the trace of spectral
projections]\label{main3} In the notation of Theorem~\ref{main1}
one has, $$ \Tr_{\Gamma}(U E(\lambda) U^*) = {\rm
rank}\,(E^0(\lambda)) \qquad {\text{for all}}\;\;\; \mu \in (0,
\mu_0). $$
\end{maincor}

The following corollary uses in addition the Rapid Decay property
(RD) for discrete groups. This property is related with the
Haagerup inequality, which estimates the convolution norm in terms
of the word lengths. Groups that are either virtually nilpotent or
word hyperbolic have property (RD). For these groups, it is also
known that every group cohomology class can be represented by a
group cocycle $c \in Z^j(\Gamma, \R)$ that is of polynomial
growth, cf. \cite{Gr}.

\begin{maincor}[Semiclassical vanishing of the higher traces of spectral
projections]\label{main2} Let $\Gamma$ be a discrete group that
has property (RD). Let $c \in Z^j(\Gamma, \R)$ ($j$ even $> 0$) be
a normalised group cocycle that is of polynomial growth, and
$\tau_c$ the induced cyclic cocycle on the twisted group algebra
$\C(\Gamma,\bar\sigma)$. Then the tensor product cocycle
$\tau_c\#\Tr$ extends continuously to $\B(\Gamma, \sigma)$, and in
the notation of Theorem~\ref{main1} one has, for all $\mu \in (0,
\mu_0)$ $$ \tau_c\#\Tr(U E(\lambda) U^*,\ldots, U E(\lambda) U^*)
= \tau_c\#\Tr(\id\otimes E^0(\lambda),\ldots, \id\otimes
E^0(\lambda))=0. $$
\end{maincor}

The method of proof of Theorem 2 is very streamlined
and almost  completely functional analytic. It uses cut-off
functions, the polar decomposition and closed image
technique. It completely replaces  a
variational  method which was used in \cite{Sh} and \cite{MS}
cited above. Though we also use the idea of model operator
\cite{Sh}, we supplement it by a direct construction
of an intertwining operator between
two spectral projections: of the original operator and the
model operator. In particular, this construction allows
to treat these two operators in a symmetric way unlike
the treatment in \cite{Sh, MS} where the proofs of the upper and
lower  estimates for the spectrum distribution functions
required separate and very different proofs.

Variational methods and cut-off functions were used
in  \cite[Proposition 5.2]{BFKM} to establish
existence of a gap near zero in the spectrum of the Witten
deformation for the periodic Laplacians on forms,
and in \cite{Bu} (see Section 4, in particular, Lemmas 4.3 and 4.10)
to prove
vanishing of the relative index term in the gluing formula
for the $\eta$-invariant in the adiabatic limit.

In \cite{HS84}, Helffer and Sj{\"o}strand used cut-off functions
and perturbation arguments based on the Riesz projection
formula (see \cite[Theorem 2.4 and Proposition 2.5]{HS84})
together with Agmon type weighted estimates to study
the tunneling effect for Schr\"odinger operators with
electric wells. These methods were extended to
magnetic Schr\"odinger operators on compact manifolds
in \cite{HS87, HS88, HS90, HS89, HSLNP, HSHaas}. It is  quite
possible that the technique developed in these papers
can also be applied to the problems studied in this
paper.

We believe, however, that our use of the
polar decomposition and closed image technique to establish
$C^*$-algebra equivalence of spectral projections is
new and may lead to further important results.
For instance, in [Kor04], the technique of this paper
was applied to prove existence of arbitrarily large
number of gaps in the spectrum of the magnetic
Schr\"odinger operators $H^h = (ih d+{\bf A})^* (ih
d+{\bf A})$ with the periodic magnetic field ${\bf
B}=d{\bf A}$ on covering spaces of compact manifolds
under some Morse type assumptions on $\bf B$ in the
semiclassical limit of strong magnetic field $h\to 0$.

We remark that as
an application of Theorems 1 and 2, one obtains a new proof of the
$L^2$-Morse inequalities, \cite{NovSh, Sh}. We also obtain
an application to the quantum Hall effect which we will now
describe.

The Kubo formula for the Hall conductance both in the usual model
of the integer quantum Hall effect on the Euclidean plane and in
the model of the fractional quantum Hall effect on the hyperbolic
plane can be naturally interpreted as a (densely defined) cyclic
2-trace $\tr_K$ on the algebra of observables ${\B}(\Gamma,
\sigma)$, \cite{Bel, CHMM, MM}. The Hall conductance
cocycle $\tr_K$ can also be shown to be given by a quadratically
bounded group cocycle. Moreover, it is well-known that $\mathbb
Z^2$ and cocompact Fuchsian groups have property (RD). Therefore
we have the following consequence of Corollary~\ref{main2}. When
$\E$ is trivial with trivial connection, the endomorphism $B$ is
zero, and the Morse type potential $V$ is a scalar valued
function, we get the magnetic Schr\"odinger operator
$H_{{\bf A},V}(\mu) = \mu^{-1} H(\mu)= \nabla^*\nabla + \mu^{-2} V$.

\begin{maincor}[Semiclassical vanishing of the Hall conductance on
low energy bands]\label{main4} Let $\M$ be either the Euclidean
plane $\R^2$ or the hyperbolic plane $\hyp$, $V$ a Morse type
potential.
Let $\lambda\in \mathbb R$ be the Fermi level,
$\lambda\not\in {\rm spec}(K)$. Let $P_\lambda =
{\chi}_{(-\infty,
\mu^{-1} \lambda]}(H_{{\bf A},V}(\mu))$ denote the spectral
projection of the magnetic Schr\"odinger operator $H_{{\bf
A},V}(\mu)$. Then there exists $\mu_0 >0$ such that for all values
of the coupling constant $\mu\in (0, \mu_0)$, the Hall conductance
vanishes, $$\sigma_\lambda = \tr_K(P_\lambda, P_\lambda,
P_\lambda) = 0.$$ That is, the low energy bands do not contribute
to the Hall conductance.
\end{maincor}

In the case of the Euclidean plane and when the magnetic field is
uniform, this result was established by a different method in
\cite{Nak+Bel}.

The physical explanation for this semiclassical vanishing theorem
for the Hall conductance is as follows. The Hall conductance for
the model operator vanishes, since it is the Hamiltonian of a
crystal with perfectly isolated atoms as mentioned earlier,
therefore there can be no current flowing through it, which
remains valid for small perturbations of the model operator.

The authors are grateful to
D. Burghelea, B. Helffer and J. Lott for helpful
comments on an earlier draft of the manuscript.

\section{Preliminaries}

Let $M$ be a compact connected Riemannian manifold, $\Gamma$ be its
fundamental group and
$\widetilde M$ be its universal cover, i.e. one has the principal bundle
$ \Gamma\to \widetilde M\overset{p}{\to} M.$ To make the paper self-contained,
we include preliminary material, some of which may not be new, cf.
\cite{Bel}, \cite{BrSu},
\cite{CHMM}, \cite{MM}, \cite{Ma}.

\subsection{Projective action, or magnetic translations}

Let $\omega$ be a closed real-valued 2-form on $M$ such that ${\bf
B} =p^* \omega$ is \emph{exact}. So ${\bf B}=d{\bf A}$ where ${\bf
A}$ is a 1-form on $\widetilde M$. We will assume without loss of
generality that ${\bf A}$ is  real-valued. Define $\nabla_{\bf
A}=d+\,i{\bf A}$.  Then $\nabla_{\bf A}$ is a Hermitian connection
on the trivial line bundle  over $\widetilde{M}$ with the
curvature $(\nabla_{\bf A})^2=i\, {\bf B}$. Suppose that $\E$ is a
Hermitian vector bundle on $M$ and $\tE$ the lift of $\E$ to the
universal cover $\widetilde M$. The connection $\nabla_{\bf A}$
defines  a projective action of $\Gamma$ on $L^2$ sections of
$\tE$ as follows.

Observe that since ${\bf B}$ is
$\Gamma$-invariant, one has
$ 0=\gamma^*{\bf B}-{\bf B}=d(\gamma^*{\bf A}-{\bf A})\quad
\forall\gamma \in\Gamma$. So $\gamma^*{\bf A}-{\bf A}$ is a closed
1-form on the
simply connected
manifold $\widetilde{M}$, therefore
\[
\gamma^*{\bf A}-{\bf A}=d\psi_\gamma,\quad\forall\gamma\in\Gamma,
\]
where $\psi_\gamma$ is a smooth
function on $\widetilde{M}$. It is defined up to an additive constant,
so we can assume in addition that  it satisfies
the following normalization condition:
\begin{itemize}
\item $\psi_\gamma(x_0)=0$ for a fixed $x_0\in\widetilde{M},\quad
\forall\gamma\in\Gamma$.
\end{itemize}
It follows that $\psi_\gamma$ is real-valued and $\psi_e(x)\equiv 0$, where
$e$ denotes the neutral element of $\Gamma$.
It is also easy to check that

\begin{itemize}
\item $\psi_\gamma(x)+\psi_{\gamma'}(\gamma
x)-\psi_{\gamma'\gamma}(x)$  is independent of
$x\in\widetilde{M},\quad \forall \gamma,\gamma'\in\Gamma$.
           \end{itemize}

Then $\sigma(\gamma,\gamma')=\exp(-i\psi_\gamma(\gamma'\cdot x_0))$  defines a
{\em multiplier} on $\Gamma$ i.e. $\sigma:\Gamma\times\Gamma\to U(1)$
satisfies
\begin{itemize}
\item ${\sigma}(\gamma,e) = {\sigma}(e,\gamma)=1,\quad\forall\
\gamma\in\Gamma$;

\item ${\sigma}(\gamma_1,\gamma_2)
{\sigma}(\gamma_1\gamma_2, \gamma_3)=
{\sigma}(\gamma_1,\gamma_2\gamma_3)
{\sigma}(\gamma_2,\gamma_3),\quad \forall \gamma_1, \gamma_2,
\gamma_3\in \Gamma$
\quad ({\em the cocycle relation}).
\end{itemize}

It follows from these relations that
$\sigma(\gamma,\gamma^{-1})=\sigma(\gamma^{-1},\gamma)$.

The complex conjugate multiplier
$\bar\sigma(\gamma,\gamma')=\exp(i\psi_\gamma(\gamma'\cdot x_0))$
also satisfies the same relations.

For $u\in L^2(\M, \tE)$ and $\gamma\in\Gamma$ define $$ \quad
U_\gamma u = (\gamma^{-1})^* u, \quad S_\gamma u =
\exp(-i\psi_\gamma)\,u. $$ Then the operators
$T_\gamma=U_\gamma\circ S_\gamma$ satisfy $$ T_e={\id},
           \quad T_{\gamma_1} T_{\gamma_2}
= {\sigma}(\gamma_1,\gamma_2) T_{\gamma_1 \gamma_2},
$$
for all $\gamma_1, \gamma_2 \in \Gamma$. In this case one says that the map
$T : \Gamma\to {\U} (L^2(\M, \tE))$, $\gamma\mapsto T_\gamma$, is a projective
$(\Gamma, {\sigma})$-unitary representation, where for any Hilbert space $\H$
           we denote by
$\U(\H)$ the group of all unitary operators in $\H$. In other words
one says that the map $\gamma\mapsto T_\gamma$ defines a
$(\Gamma,\sigma)$-{\em action}
in~$\H$.

It is also easy to check that the adjoint operator to $T_\gamma$
in $L^2(\M, \tE)$ (with respect to a smooth $\Gamma$-invariant
measure on $\M$ and a $\Gamma$-invariant Hermitian structure on
$\tE$) is $$
T_\gamma^*=\bar\sigma(\gamma,\gamma^{-1})T_{\gamma^{-1}}. $$

The operators $T_\gamma$ are also called
{\em magnetic translations}.

\subsection{Twisted group algebras}
Denote by $\ell^2(\Gamma)$ the standard Hilbert space of
complex-valued $L^2$-functions on the discrete group $\Gamma$. We
will use a left $(\Gamma, \bar\sigma)$-action on $\ell^2(\Gamma)$
(or, equivalently, a $(\Gamma, \bar\sigma)$-unitary representation
in $\ell^2(\Gamma)$) which is given explicitly by $$ T_\gamma^L
f(\gamma') = f(\gamma^{-1}\gamma') \bar\sigma(\gamma,
\gamma^{-1}\gamma'), \qquad \gamma, \gamma' \in \Gamma. $$ It is
easy to see that this is indeed a $(\Gamma,\bar\sigma)$-action,
i.e. $$ T_e^L={\id} \quad {\rm and} \quad T_{\gamma_1}^L
T_{\gamma_2}^L= \bar\sigma(\gamma_1,\gamma_2)T_{\gamma_1
\gamma_2}^L, \quad \forall \gamma_1,\gamma_2\in\Gamma. $$ Also $$
(T_\gamma^L)^*=\sigma(\gamma, \gamma^{-1})T_{\gamma^{-1}}^L. $$

Let
$$
{\mathcal A}^R(\Gamma, \sigma) =
\Big\{ A \in \B(\ell^2(\Gamma)): [T_\gamma^L, A] = 0, \quad \forall \gamma \in
\Gamma\Big\}
$$
be the commutant of the left $(\Gamma, \bar\sigma)$-action on
$\ell^2(\Gamma)$.  Here by $\B(\H)$ we denote the algebra of all
bounded linear operators
in a Hilbert space $\H$.
By the general theory, ${\mathcal A}^R(\Gamma, \sigma)$ is a von
Neumann algebra
and is known as the {\em (right) twisted group von Neumann algebra}.
It can also be realized as follows. Let us define the  following operators
in $\ell^2(\Gamma)$: $$ T_\gamma^R f(\gamma') = f(\gamma'\gamma)
\sigma(\gamma', \gamma), \qquad  \gamma, \gamma' \in \Gamma. $$ It
is easy to check that they form a right $(\Gamma,\sigma)$-action
in $\ell^2(\Gamma)$ i.e. $$ T_e^R={\rm id} \quad {\rm and} \quad
T_{\gamma_1}^R  T_{\gamma_2}^R=
\sigma(\gamma_1,\gamma_2)T_{\gamma_1 \gamma_2}^R, \quad \forall
\gamma_1,\gamma_2\in\Gamma, $$ and also $$
(T_{\gamma}^R)^*=\bar\sigma(\gamma,\gamma^{-1}) T_{\gamma^{-1}}^R.
$$ This action commutes with the left $(\Gamma,\bar\sigma)$-action
defined above i.e. $$ T_\gamma^L  T_{\gamma'}^R= T_{\gamma'}^R
T_\gamma^L, \quad \forall \gamma,\gamma'\in\Gamma. $$ It can be
shown that  the von Neumann algebra $\A^R(\Gamma,\sigma)$ is
generated by the operators $\{T_{\gamma}^R\}_{\gamma\in\Gamma}$
(see e.g.\ a similar argument in \cite{Sh2}).

Similarly we can introduce  a von Neumann algebra
$$
{\mathcal A}^L(\Gamma, \bar\sigma) =
\Big\{ A \in \B(\ell^2(\Gamma)): [T_\gamma^R, A] = 0, \quad \forall \gamma \in
\Gamma\Big\}.
$$
We will refer to it as {\em (left) twisted group von Neumann algebra}.
It is generated by the operators
$\{T_{\gamma}^L\}_{\gamma\in\Gamma}$, and it is the commutant
of ${\mathcal A}^R(\Gamma, \sigma)$.

Let us define a twisted group algebra $\C(\Gamma,\sigma)$ which consists of
complex valued functions with finite support on $\Gamma$ and with the
twisted convolution
operation
$$
(f*g)(\gamma)=\sum_{\gamma_1,\gamma_2:\gamma_1\gamma_2=
\gamma}f(\gamma_1)g(\gamma_2)\sigma(\gamma_1,\gamma_2).
$$
The basis of $\C(\Gamma,\sigma)$ as a vector space is formed by
$\delta$-functions
$\{\delta_\gamma\}_{\gamma\in\Gamma}$, $\delta_\gamma(\gamma')=1$ if
$\gamma=\gamma'$
and $0$ otherwise. We have
$$
\delta_{\gamma_1} *
\delta_{\gamma_2}=\sigma(\gamma_1,\gamma_2)\delta_{\gamma_1\gamma_2}.
$$
Associativity of this multiplication is equivalent to the cocycle condition.

Note also that the $\delta$-functions $\{\delta_\gamma\}_{\gamma\in\Gamma}$
form an orthonormal basis in $\ell^2(\Gamma)$. It is easy to check that
$$
T_\gamma^L\delta_{\gamma'}=\bar\sigma(\gamma,\gamma')\delta_{\gamma\gamma'},
\quad
T_\gamma^R\delta_{\gamma'}=
\sigma(\gamma'\gamma^{-1},\gamma)\delta_{\gamma'\gamma^{-1}}.
$$

It is clear that the correspondences $\delta_\gamma\mapsto T^L_{\gamma}$
and $\delta_\gamma\mapsto  T^R_{\gamma}$ define representations
of $\C(\Gamma,\bar\sigma)$ and $\C(\Gamma,\sigma)$ respectively. In both cases
the weak closure of the image of the twisted group algebra coincides
with the corresponding
von Neumann algebra ($\A^L(\Gamma,\bar\sigma)$ and
$\A^R(\Gamma,\sigma)$ respectively). The
corresponding norm closures are so called {\em reduced twisted group}
$C^*$-{\em algebras} which are denoted $C^*_r(\Gamma,\bar\sigma)$ and
$C^*_r(\Gamma,\sigma)$
respectively.

The von Neumann algebras $\A^L(\Gamma,\bar\sigma)$ and
$\A^R(\Gamma,\sigma)$ can be described
in terms of the matrix elements. For any $A\in\B(\ell^2(\Gamma))$ denote
$A_{\alpha,\beta}=(A\delta_\beta,\delta_\alpha)$ (which is a matrix
element of $A$).
Then repeating standard arguments (given in a similar situation e.g.\
in \cite{Sh2}) we
can prove that for any $A\in\B(\ell^2(\Gamma))$ the inclusion $A\in
\A^R(\Gamma,\sigma)$
is equivalent to the relations
$$
A_{\gamma x,\gamma y}=\bar\sigma(\gamma,x)\sigma(\gamma,y)A_{x,y}\;,\quad
\forall x,y,\gamma\in\Gamma.
$$
In particular, we have for any $A\in \A^R(\Gamma,\sigma)$
$$
A_{\gamma x,\gamma x}=A_{x,x}\;, \quad \forall x,\gamma\in\Gamma.
$$
Similarly, for any $A\in\B(\ell^2(\Gamma))$ the inclusion $A\in
\A^L(\Gamma,\bar\sigma)$
is equivalent to the relations
$$
A_{x\gamma,y\gamma}=\bar\sigma(x,\gamma)\sigma(y,\gamma)A_{x,y}\;,\quad
\forall x,y,\gamma\in\Gamma.
$$
In particular, we have
$$
A_{x\gamma ,x\gamma }=A_{x,x}\;, \quad \forall x,\gamma\in\Gamma,
$$
for any $A\in \A^L(\Gamma,\bar\sigma)$.

A finite {\em von Neumann trace}
$\tr_{\Gamma,\bar\sigma}:\A^L(\Gamma,\bar\sigma)\to\C$
is defined by the formula
$$
\tr_{\Gamma,\bar\sigma} A=(A\delta_e,\delta_e).
$$
We can also write $\tr_{\Gamma,\bar\sigma} A=A_{\gamma,\gamma}=
\left(A\delta_\gamma,\delta_\gamma\right)$  for any $\gamma\in
\Gamma$ because the right hand
side does not depend of
$\gamma$.

A finite von Neumann trace $\tr_{\Gamma,\sigma}:\A^R(\Gamma,\sigma)\to\C$
is defined by the same formula, so we will denote by $\tr_\Gamma$ any
of these traces.

Let $\mathcal H$ denote an infinite dimensional complex Hilbert
space. Then the Hilbert tensor product $\ell^2(\Gamma)\otimes
\mathcal H$ is both $(\Gamma, \bar\sigma)$-module and $(\Gamma,
\sigma)$-module under the actions $\gamma \mapsto T_\gamma^L
\otimes \id$ and $\gamma \mapsto  T_\gamma^R \otimes \id$
respectively. Let ${\mathcal A}^L_{\mathcal H}(\Gamma,\bar\sigma)$
and ${\mathcal A}^R_{\mathcal H}(\Gamma,\sigma)$ denote the von
Neumann algebras in $\ell^2(\Gamma)\otimes \mathcal H$ which are
commutants of the $(\Gamma, \sigma)$- and $(\Gamma,
\bar\sigma)$-actions respectively. Clearly ${\mathcal
A}^L_{\mathcal H}(\Gamma,\bar\sigma) \cong {\mathcal
A}^L(\Gamma,\bar\sigma) \otimes \B(\mathcal H)$ and ${\mathcal
A}^R_{\mathcal H}(\Gamma,\sigma) \cong {\mathcal
A}^R(\Gamma,\sigma) \otimes \B(\mathcal H)$ in the usual sense of
von Neumann algebra tensor products. Moreover, we have the
following

\begin{lemma}
Any operator $A\in {\mathcal A}^L_{{\mathcal
H}}(\Gamma,\bar\sigma)$ can be represented as
\[
A=\sum_{\gamma\in\Gamma}T^L_\gamma\otimes A(\gamma),
\]
where $A(\gamma)\in \B({\mathcal H})$, and the series in the
right-hand side of this identity converges in the strong operator
topology.
\end{lemma}

\begin{proof}
Let $A\in {\mathcal A}^L_{\mathcal H}(\Gamma,\bar\sigma)$. Define
a bounded operator $A(\gamma)$ in ${\mathcal H}$ by the formula
\[
A(\delta_e\otimes v)=\sum_{\gamma\in\Gamma}\delta_\gamma\otimes
A(\gamma)v, \quad v\in {\mathcal H}.
\]
Take any $x=\sum_{\gamma\in\Gamma}\delta_\gamma\otimes x_\gamma
\in \ell^2(\Gamma)\otimes {\mathcal H}$. Then we have
\begin{eqnarray*}
Ax&=& \sum_{\gamma_1} A(\delta_{\gamma_1}\otimes x_{\gamma_1})\\
&=& \sum_{\gamma_1} \sigma(\gamma_1,\gamma_1^{-1})^{-1}
A(T^R_{\gamma_1^{-1}}\otimes \id)(\delta_{e}\otimes
x_{\gamma_1})\\ &=& \sum_{\gamma_1}
\sigma(\gamma_1,\gamma_1^{-1})^{-1} (T^R_{\gamma_1^{-1}}\otimes
\id)A(\delta_{e}\otimes x_{\gamma_1})\\ &=& \sum_{\gamma_1}
\sigma(\gamma_1,\gamma_1^{-1})^{-1} (T^R_{\gamma_1^{-1}}\otimes
\id) \sum_{\gamma_2}\delta_{\gamma_2}\otimes
A(\gamma_2)x_{\gamma_1}\\ &=& \sum_{\gamma_1,\gamma_2}
\sigma(\gamma_1,\gamma_1^{-1})^{-1}
\sigma(\gamma_2\gamma_1,\gamma_1^{-1})
\delta_{\gamma_2\gamma_1}\otimes A(\gamma_2)x_{\gamma_1}.
\end{eqnarray*}
By the cocycle identity, we have $\sigma(\gamma_2,\gamma_1)
\sigma(\gamma_2\gamma_1,\gamma^{-1}_1) =\sigma(\gamma_2,e)
\sigma(\gamma_1,\gamma_1^{-1})$, that implies
\[
\sigma(\gamma_1,\gamma_1^{-1})^{-1}\sigma(\gamma_2\gamma_1,\gamma_1^{-1})=
\bar{\sigma}(\gamma_2,\gamma_1)
\]
and, finally, gives
\[
Ax=\sum_{\gamma_1,\gamma_2}\bar{\sigma}(\gamma_2,\gamma_1)
\delta_{\gamma_2\gamma_1}\otimes A(\gamma_2)x_{\gamma_1}=
\sum_{\gamma_1,\gamma_2} T^L_{\gamma_2}\delta_{\gamma_1}\otimes
A(\gamma_2)x_{\gamma_1},
\]
that completes the proof.
\end{proof}

Let us note the following elementary lemma.

\begin{lemma}\label{C}
Let $A\in {\mathcal A}^L_{\mathcal H}(\Gamma,\bar\sigma)$,
$A=\sum_{\gamma\in\Gamma}T^L_\gamma\otimes A(\gamma)$, where
$A(\gamma)\in \B({\mathcal H})$. Then
\[
\sup_{\ga\in\Gamma}\|A(\gamma)\|\le \|A\| \le \sum_{\gamma\in
\Gamma} \|A(\gamma)\|,
\]
where the righthand side of the inequality is not necessarily
finite.
\end{lemma}

Define the semifinite tensor product trace $\Tr_{\Gamma} =
\tr_{\Gamma} \otimes \Tr$ on each of the algebras ${\mathcal
A}^L_{\mathcal H}(\Gamma,\bar\sigma)$ and ${\mathcal
A}^R_{\mathcal H}(\Gamma,\sigma)$. Here $\Tr$ denotes the standard
(semi-finite) trace on $\B(\H)$.

The $C^*$-tensor product $C^*_r(\Gamma,
\bar\sigma)\otimes\K(\mathcal H)$ is the norm closure of the
algebraic tensor product $C^*_r(\Gamma,
\bar\sigma)\odot\K(\mathcal H)\subset {\mathcal A}^L_{\mathcal
H}(\Gamma,\bar\sigma) \cong {\mathcal A}^L(\Gamma,\bar\sigma)
\otimes \B(\mathcal H)$ in $\B(\ell^2(\Gamma)\otimes\mathcal H)$
One can give the following sufficient conditions for an operator
$A\in {\mathcal A}^L_{\mathcal H}(\Gamma,\bar\sigma)$ to belong to
the algebra $C^*_r(\Gamma,\bar\sigma) \otimes\mathcal K(\mathcal
H)$.

\begin{lemma}\label{D1}
If $A\in {\mathcal A}^L_{\mathcal H}(\Gamma,\bar\sigma)$,
$A=\sum_{\gamma\in\Gamma}T^L_\gamma\otimes A(\gamma)$ is such that
$A(\gamma) \in \K(\mathcal H)$ and also satisfies
$
\;\displaystyle\sum_{\gamma} \|A(\gamma)\|<\infty, $ then $A \in
C^*_r(\Gamma, \bar\sigma) \otimes\K(\mathcal H)$ and $\|A\| \le
\;\displaystyle\sum_{\gamma} \| A(\gamma)\|$.
\end{lemma}

\begin{proof}
Let $K_1 \subset K_2 \subset \cdots $ be a sequence of finite
subsets of $\Gamma$ which is an exhaustion of $\Gamma$, i.e.
$\bigcup_{j\ge 1} K_j = \Gamma$. For all $j \in \mathbb N$, define
$A_j \in {\mathcal A}^L_{\mathcal H}(\Gamma,\bar\sigma)$ as
$A_j=\sum_{\gamma\in\Gamma}T^L_\gamma\otimes A_j(\gamma)$, where
$$ { {A_j}} (\gamma) = \left\{\begin{array}{l} {{A}} (\gamma) \;\;
{\rm if} \;\; \gamma\in K_j;\\[7pt]
         0 \;\;  {\rm otherwise.}
\end{array}\right.
$$ Then in fact $A_j\in \C(\Gamma,\bar\sigma)\otimes {\mathcal
K(\mathcal H)}$ by definition. Using Lemma~\ref{C}, we have $$
\begin{array}{lcl}
\|A-A_j\| & \le & \displaystyle\sum_{\gamma} \|({A -
A_j})(\gamma)\|\\[+7pt] & = & \displaystyle\sum_{\gamma}
\|{A}(\gamma) - {A_j}(\gamma)\|\\[+7pt] & = &
\displaystyle\sum_{\gamma\in \Gamma\setminus K_j} \|A(\gamma)\|.
\end{array}
$$ By hypothesis, $\displaystyle\sum_{\gamma}
\|A(\gamma)\|<\infty$, therefore $\displaystyle\sum_{\gamma\in
\Gamma\setminus K_j} \|A(\gamma)\| \to 0$ as $j\to\infty$, since
$K_j$ is an increasing exhaustion of $\Gamma$. This proves that $A
\in  C^*_r(\Gamma, \bar\sigma)\otimes\K(\mathcal H)$.
\end{proof}

\subsection{Projectively invariant elliptic operators}

As before, let $M$ be a closed Riemannian manifold and $\widetilde
M$ be its universal cover. Let $\E$ be a Hermitian vector bundle
on $M$ and $\tE$ the lift of $\E$ to the universal cover
$\widetilde M$. Let $\widetilde \nabla^\E$ denote a
$\Gamma$-invariant Hermitian connection on $\tE$. Then consider
the Hermitian connection $\nabla = \widetilde\nabla^\E\otimes \id
+ \id \otimes \nabla_{\bf A}$ on $\tE\otimes {\mathcal L} = \tE$,
where ${\mathcal L}$ is the trivial line bundle on $\M$,
$\nabla_{\bf A}=d+i{\bf A}$.

Using the Riemannian metric on $\widetilde M$ and the Hermitian
metric on $\tE$, consider the elliptic self-adjoint differential
operator given by
\begin{equation}\label{ham2}
H(\mu) = \mu \nabla^*\nabla + B  + \mu^{-1} V,
\end{equation}
where $B$, $V$ are $\Gamma$-invariant self-adjoint endomorphisms of
the  bundle $\tE$, $\mu$ is the coupling constant
and where $V$ satisfies
in addition the Morse type condition.
Then $H(\mu)$ acts on $L^2(\M, \tE)$ and is a self-adjoint second order
elliptic differential
operator. It commutes with the  magnetic  translations $T_\gamma$ (for all
$\gamma\in\Gamma$), i.e.  with the $(\Gamma,{\sigma})$-action which
was defined above. To see this note first that the operators
$U_\gamma=(\gamma^{-1})^*$ and $S_\gamma$ (the multiplication by
$\exp(-i\psi_\gamma)$)
are defined not only on sections of $\tE$ but also on $\tE$-valued 1-forms
(and actually on $\tE$-valued $p$-forms for any $p\ge 0$) on $\widetilde M$.
Hence the magnetic translations $T_\gamma$ are well defined on
$\tE$-valued forms as well.
The operators $T_\gamma$ are obviously unitary on the $L^2$ spaces of
$\tE$-valued forms,
where the $L^2$ structure is defined by the fixed $\Gamma$-invariant metric on
$\widetilde M$ and the fixed $\Gamma$-invariant Hermitian metric on $\tE$.
An easy calculation shows that $T_\gamma\nabla=\nabla T_\gamma$
on sections of $\tE$. By taking adjoint operators we obtain
$T_\gamma \nabla^*=\nabla^* T_\gamma$ on $\tE$-valued 1-forms.
Therefore $T_\gamma  \nabla^*\nabla= \nabla^*\nabla T_\gamma$
on sections of $\tE$. Since obviously $T_\gamma B=B T_\gamma$ and
$T_\gamma V=V T_\gamma$, we see that $H(\mu)$
commutes with $T_\gamma$ for all $\gamma$.

Since $H(\mu)$ commutes with the $(\Gamma,{\sigma})$-action, it
follows by the spectral mapping theorem that  the spectral
projections of $H(\mu)$, $E(\lambda) =
\chi_{(-\infty,\lambda]}(H(\mu))$ are bounded operators on
$L^2(\M, \tE)$ that also commute with the
$(\Gamma,{\sigma})$-action i.e. $T_\gamma E(\lambda)= E(\lambda)
T_\gamma,\quad\forall\ \gamma\in\Gamma$. The {\em commutant} of
the $(\Gamma, {\sigma})$-action is a von Neumann algebra
\[
{\mathcal U}_{\mathfrak{H}}(\Gamma,\bar\sigma) = \left\{Q\in
\B(\mathfrak{H}) : T_\gamma Q = Q T_\gamma,\quad\forall\
\gamma\in\Gamma\right\},
\] where $\mathfrak{H} = L^2(\M, \tE)$.
To characterize the Schwartz kernels $k_Q(x,y)$ of the operators $Q\in
{\mathcal U}_{\mathfrak{H}}(\Gamma,\bar\sigma)$ note that the relation
$T_\gamma Q=QT_\gamma$ can be rewritten in the form
\begin{equation}\label{kernel}
e^{i\psi_\gamma(x)} k_Q(\gamma x,\gamma y) e^{-i\psi_\gamma(y)} =
k_Q(x,y), \quad \forall x,y\in\M\quad \forall \gamma\in\Gamma,
\end{equation} where we have identified the fibre $\tE_x$ with the
fibre at $\tE_{\gamma x}$ via the isomorphism induced by $\gamma$.
So $Q\in {\mathcal U}_{\mathfrak{H}}(\Gamma,\bar\sigma)$ if and
only if $Q\in \B(\mathfrak{H})$ and (\ref{kernel}) holds. In
particular, in this case $k_Q(x,x)$ is $\Gamma$-invariant. For the
spectral projections of $H(\mu)$ we also have $E(\lambda)\in
{\mathcal U}_{\mathfrak{H}}(\Gamma,\bar\sigma)$, so the
corresponding Schwartz kernels also satisfy (\ref{kernel}). Note
that by elliptic regularity, the Schwartz kernels of $E(\lambda)$
are smooth.

To define a natural trace on ${\mathcal U}_{\mathfrak{H}}(\Gamma,\bar\sigma)$
we will construct an isomorphism of this algebra
with the von Neumann algebra $\A^L_\H(\Gamma,\bar\sigma)$, where
$\H=L^2(\F, \tE|_\F)$ and
       $\mathcal{F}$ is a fundamental domain for the $\Gamma$-action on
$\M$. By choosing a connected fundamental domain $\mathcal{F}$ for
the action of
$\Gamma$ on $\M$, we can
define a $(\Gamma, \sigma)$-equivariant isometry
\begin{equation}\label{isom}
{\bf U} : L^2(\M, \tE)\cong \ell^2(\Gamma)\otimes L^2(\mathcal{F},  \tE|_\F)
\end{equation}
as follows. Let $i: \mathcal{F} \to \M$ denote the inclusion map.
Define
$$
{\bf U} (\phi) = \sum_{\gamma\in \Gamma}\delta_\gamma\otimes
i^*(T_\gamma \phi),
\qquad  \phi \in L^2(\M, \tE).
$$

\begin{lemma} The map
$\; {\bf U} :  L^2(\M, \tE)\to  \ell^2(\Gamma)\otimes
L^2(\mathcal{F},  \tE|_\F)\;$ defined
above in (\ref{isom}) is a
           $(\Gamma, \sigma)$-equivariant unitary operator, where the
$(\Gamma, \sigma)$-action on $\ell^2(\Gamma)\otimes
L^2(\mathcal{F}, \tE|_\F)$ is given by the operators
$T^R_\gamma\otimes \id$.
\end{lemma}

\begin{proof} Given $\phi \in L^2(\M, \tE)$, we compute

\begin{align*}
{\bf U}(T_\gamma\phi)
&=\sum_{\gamma'\in\Gamma}\delta_{\gamma'}\otimes
i^*(T_{\gamma'}T_\gamma\phi)
=\sum_{\gamma'\in\Gamma}\sigma(\gamma',\gamma)\delta_{\gamma'}\otimes
i^*(T_{\gamma'\gamma}\phi)\\
&=\sum_{\gamma'\in\Gamma}\sigma(\gamma'\gamma^{-1},\gamma)\delta_{\gamma'\gamma^{-1}}
\otimes i^*(T_{\gamma'}\phi) =(T^R_\gamma\otimes \id){\bf U}\phi,
\end{align*}
which proves that $\bf U$ is a $(\Gamma,\sigma)$-equivariant map.
It is straightforward to check that the operator $\bf U$ is unitary.
\end{proof}

Since $\mathcal{U}_{\mathfrak{H}}(\Gamma,\bar\sigma)$ is the
commutant of $\{T_\gamma\}_{\gamma\in\Gamma}$, and
$\A^L_\H(\Gamma,\bar\sigma)$ is the commutant of
$\{T^R_\gamma\otimes \id\}_{\gamma\in\Gamma}$, we see that $\bf U$
induces an isomorphism of von Neumann algebras
$\mathcal{U}_{\mathfrak{H}}(\Gamma,\bar\sigma)$ and
$\A^L_\H(\Gamma,\bar\sigma)$. Therefore we can transfer the trace
$\Tr_\Gamma$ from $\A^L_\H(\Gamma,\bar\sigma)$ to
$\mathcal{U}_{\mathfrak{H}}(\Gamma,\bar\sigma)$. The result will
be a semifinite $\Gamma$-trace on
$\mathcal{U}_{\mathfrak{H}}(\Gamma,\bar\sigma)$ which we will
still denote $\Tr_\Gamma$.

Note that the trace $\Tr_\Gamma$ is faithful. This follows by well 
known arguments
from the theory of von Neumann algebras, which are reproduced
for instance, in \cite{Ta} on pages 316-317,
in the proof of Proposition V.2.14.
The result in \cite{Ta} is directly applicable to the case of the trivial
multiplier only, but the arguments easily work in the general case.
In the case of non-trivial multiplier the trace was briefly
considered by Br\"uning and Sunada \cite{BrSu},
who also observed the faithfulness of the trace.

It is easy to check that for any $Q\in
\mathcal{U}_{\mathfrak{H}}(\Gamma,\bar\sigma)$ with a finite
$\Gamma$-trace and a continuous Schwartz kernel $k_Q$ we have $$
\Tr_\Gamma Q=\int_\F \tr k_Q(x,x)dx $$ where $dx$ denotes the
$\Gamma$-invariant measure and $\tr$ the pointwise trace. An
important particular case is a spectral projection $E(\lambda)$ of
the elliptic self-adjoint operator $H(\mu)$, which has a smooth
Schwartz kernel and so a finite $\Gamma$-trace. Therefore we can
define a {\em spectral density function} $$
N_\Gamma(\lambda;H)=\Tr_\Gamma E(\lambda), $$ which is finite for
all $\lambda\in \R$. It is easy to see that $\lambda\mapsto
N_\Gamma(\lambda;H(\mu))$ is a non-decreasing function, and the
spectrum of $H(\mu)$ can be reconstructed  as the set of its
points of growth, i.e. $$ \spec (H(\mu))=\{\lambda\in\R:
N_\Gamma(\lambda+\eps;H(\mu))-N_\Gamma(\lambda-\eps;H(\mu))>0, \
\forall \eps>0\}. $$

\section{Refined semiclassical approximation principle and existence
of spectral gaps} \label{equivalence} The main goal of this
section is to prove Theorem~\ref{main0} and the first part of
Theorem~\ref{main1}. We start with an abstract operator-theoretic
setting, where similar results can be stated. Then these results
are applied to projectively invariant elliptic operators with
invariant Morse type potentials on covering spaces of compact
manifolds.

\subsection{General results on equivalence of projections and existence of
spectral gaps}\label{s:abstract-equiv} Let $\mathfrak A$ be a
$C^*$-algebra, $\chh$ a Hilbert space equipped with a faithful
$\ast$-representation of ${\mathfrak A}$, $\pi: {\mathfrak A}\to
\ck(\chh)$. For simplicity of notation, we will often identify the
algebra ${\mathfrak A}$ with its image $\pi({\mathfrak A})$.

Consider Hilbert spaces ${\mathcal H}_1$ and ${\mathcal H}_2$
equipped with inner products $(\cdot,\cdot)_1$ and
$(\cdot,\cdot)_2$. 
Assume that there are given unitary operators
${\mathcal V}_1 : 
{\mathcal H}_1\to {\mathcal H}$ and ${\mathcal
V}_2 : {\mathcal 
H}_2\to {\mathcal H}$. Using the unitary
isomorphisms ${\mathcal 
V}_1$ and ${\mathcal V}_2$, we get
representations $\pi_1$ and 
$\pi_2$ of ${\mathfrak A}$ in
${\mathcal H}_1$ and ${\mathcal H}_2$ 
accordingly,
$\pi_l(a)={\mathcal V}^{-1}_l\circ\pi(a)\circ {\mathcal 
V}_l,
l=1,2, a\in {\mathfrak A}$.

Consider (unbounded) self-adjoint 
operators $A_1$ in ${\mathcal
H}_1$ and $A_2$ in ${\mathcal H}_2$ 
with the domains $\Dom (A_1)$
and $\Dom (A_2)$ respectively. We will 
assume that
\begin{itemize}
\item the operators $A_1$ and $A_2$ are 
semi-bounded from below:
\begin{gather}\label{e:5}
(A_1u,u)_1\geq 
\lambda_{01}\|u\|_1^2,\quad u\in \Dom (A_1), \\
\label{e:6} 
(A_2u,u)_2\geq \lambda_{02}\|u\|_2^2,\quad u\in \Dom
(A_2),
\end{gather}
with 
some $\lambda_{01},\lambda_{02}\leq 0$;
\item for any $t>0$, the 
operators $e^{-tA_l}, l=1,2,$ belong to
$\pi_l({\mathfrak 
A})$.
\end{itemize}

Let ${\mathcal H}_0$ be a Hilbert space, 
equipped with injective
bounded linear maps $i_1:{\mathcal H}_0\to 
{\mathcal H}_1$ and
$i_2:{\mathcal H}_0\to {\mathcal H}_2$. Assume 
that there are
given bounded linear maps $p_1:{\mathcal H}_1\to 
{\mathcal H}_0$
and $p_2:{\mathcal H}_2\to {\mathcal H}_0$ such that 
$p_1\circ
i_1=\id_{{\mathcal H}_0}$ and $p_2\circ i_2=\id_{{\mathcal H}_0}$.
The whole picture can be represented by the following diagram
(note that this diagram is not commutative).

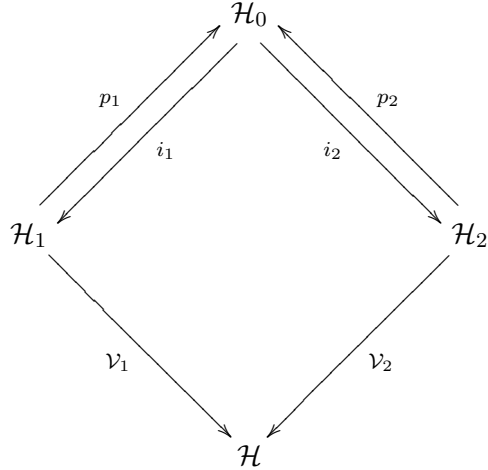
\begin{equation*}
\xymatrix @=8pc @ur { \cH_1 \ar@<1ex>[r]^{p_1}\ar[d]_{{\mathcal
V}_1} & \cH_0\ar@<-1ex>[d]_{i_2} \ar@<1ex>[l]^{i_1} \\ \cH & \cH_2
\ar@<-1ex>[u]_{p_2} \ar[l]^{{\mathcal V}_2}}
\end{equation*}

Consider a self-adjoint bounded operator $J$ in ${\mathcal H}_0$.
We assume that
\begin{itemize}
\item the operator ${\mathcal V}_2i_2Jp_1{\mathcal
V}^{-1}_1$ belongs to the von Neumann algebra $\pi({\mathfrak
A})''$;
\item $(i_2Jp_1)^*=i_1Jp_2$;
\item for any $a\in {\mathfrak A}$, the operator 
$\pi(a){\mathcal
V}_2(i_2Jp_1){\mathcal V}^{-1}_1$ belongs to 
$\pi({\mathfrak A})$.
\end{itemize}

Since the operators 
$i_l:{\mathcal H}_0\to {\mathcal H}_l, l=1,2,$
are bounded and have 
bounded left-inverse operators $p_l$, they
are topological 
monomorphisms, i.e. they have closed image and the
maps 
$i_l:{\mathcal H}_0\to {\rm Im}\, i_l$ are topological
isomorphisms. 
Therefore, we can assume that the 
estimate
\begin{equation}\label{e:rho}
\rho^{-1}\|i_2Ju\|_2\leq 
\|i_1Ju\|_1\leq \rho\|i_2Ju\|_2, \quad
u\in {\mathcal 
H}_0,
\end{equation}
holds with some $\rho>1$ (depending on 
$J$).

Define the bounded operators $J_l$ in ${\mathcal H}_l, l=1,2,$ 
by
the formula $J_l=i_lJp_l$. We assume that
\begin{itemize}
\item 
the operator $J_l, l=1,2,$ maps the domain of $A_l$ to itself;
\item 
$J_l$ is self-adjoint, and $0\leq J_l\leq \id_{{\mathcal H}_l}, 
l=1,2$;
\item for $u\in {\mathcal H}_0$, $i_1Ju\in \Dom(A_1)$ iff 
$i_2Ju\in \Dom(A_2)$.
\end{itemize}
Denote $D=\{u\in {\mathcal H}_0 : 
i_1Ju\in \Dom(A_1)\}=\{u\in
{\mathcal H}_0 : i_2Ju\in 
\Dom(A_2)\}.$

Introduce a self-adjoint positive bounded linear 
operator $J'_l$
in ${\mathcal H}_l$ by the 
formula
$J_l^2+{J'_l}{}^2=\id_{{\mathcal H}_l}$. We assume 
that
\begin{itemize}
\item the operator $J'_l, l=1,2,$ maps the 
domain of $A_l$ to itself;
\item the operators $[J_l,[J_l,A_l]]$ and 
$[J'_l,[J'_l,A_l]]$ extend to bounded
operators in ${\mathcal H}_l$, 
and
\begin{equation}\label{e:15}
\max (\|[J_l,[J_l,A_l]]\|_l,\, 
\|[J'_l,[J'_l,A_l]]\|_l)\leq
\gamma_l, \quad 
l=1,2.
\end{equation}
\end{itemize}

Finally, we assume 
that
\begin{equation}\label{e:14}
(A_lJ'_lu,J'_lu)_l\geq \alpha_l 
\|J'_lu\|_l^2,\quad u\in
\Dom(A_l), \quad l=1,2,
\end{equation}
for 
some $\alpha_l>0$, 
and
\begin{align}\label{e:16}
(A_2i_2Ju,i_2Ju)_2\leq 
\beta_1
(A_1i_1Ju,i_1Ju)_1+\varepsilon_1\|i_1Ju\|_1^2, \quad u\in D, 
\\
\label{e:A1A2} (A_1i_1Ju,i_1Ju)_1\leq 
\beta_2
(A_2i_2Ju,i_2Ju)_2+\varepsilon_2\|i_2Ju\|_2^2, \quad u\in 
D,
\end{align}
for some $\beta_1,\beta_2\geq 1$ and 
$\varepsilon_1,
\varepsilon_2>0$.

Denote by $E_l(\lambda), l=1,2$, 
the spectral projection of the
operator $A_l$, corresponding to the 
semi-axis
$(-\infty,\lambda]$. We assume that there exists a 
faithful,
normal, semi-finite trace $\tau$ on $\pi({\mathfrak A})''$ 
such
that, for any $t>0$, the operators 
${\mathcal
V}_le^{-tA_l}{\mathcal V}^{-1}_l, l=1,2,$ belong 
to
$\pi({\mathfrak A})$ and have finite trace. By standard 
arguments,
it follows that ${\mathcal V}_lE_l(\lambda){\mathcal 
V}^{-1}_l\in
\pi ({\mathfrak A})''$, and 
$\tau({\mathcal
V}_lE_l(\lambda){\mathcal V}^{-1}_l)<\infty$ for any 
$\lambda,
l=1,2$.

\begin{thm}\label{t:equivalence}
Under current 
assumptions, let $b_1>a_1$ 
and
\begin{align}\label{e:a2}
a_2&=\rho\left[ \beta_1 
\left(a_1+\gamma_1+
\frac{(a_1+\gamma_1-\lambda_{01})^2}{\alpha_1-a_1-\gamma_1}\right)+
\varepsilon_1\right],\\ \label{e:b2}
b_2&=\frac{\beta_2^{-1}(b_1\rho^{-1}-\varepsilon_2)(\alpha_2-\gamma_2)
-\alpha_2\gamma_2+2\lambda_{02}\gamma_2-\lambda^2_{02}}
{\alpha_2-2\lambda_{02}+\beta_2^{-1}(b_1\rho^{-1}-\varepsilon_2)}.
\end{align}
Suppose that $\alpha_1>a_1+\gamma_1$, $\alpha_2>b_2+\gamma_2$ and
$b_2>a_2$. If the interval $(a_1,b_1)$ does not intersect with the
spectrum of $A_1$, then:

(1) the interval $(a_2,b_2)$ does not intersect with the spectrum
of $A_2$;

(2) for any $\lambda_1\in (a_1,b_1)$ and $\lambda_2\in (a_2,b_2)$,
the projections ${\mathcal V}_1E_1(\lambda_1){\mathcal V}^{-1}_1$
and ${\mathcal V}_2E_2(\lambda_2){\mathcal V}^{-1}_2$ belong to
${\mathfrak A}$ and are Murray-von Neumann equivalent in
${\mathfrak A}$.
\end{thm}

\newtheorem{rem}[thm]{Remark}

\begin{rem}\label{b1b2}
Since $\rho>1, \beta_1\geq 1, \gamma_1>0$ and $\varepsilon_1>0$,
we, clearly, have $a_2>a_1$. The formula (\ref{e:b2}) is
equivalent to the formula
\[
b_1=\rho\left[ \beta_2 
\left(b_2+\gamma_2+
\frac{(b_2+\gamma_2-\lambda_{02})^2}{\alpha_2-b_2-\gamma_2}\right)+
\varepsilon_2\right],
\]
which 
is obtained from (\ref{e:a2}), if we replace $\alpha_1,
\beta_1, \gamma_1, \varepsilon_1, \lambda_{01}$ by $\alpha_2,
\beta_2, \gamma_2, \varepsilon_2, \lambda_{02}$ accordingly and
$a_1$ and $a_2$ by $b_2$ and $b_1$ accordingly. In particular,
this implies that $b_1>b_2$.
\end{rem}

\subsection{Localization theorem for spectral projections}
The goal of this Section is to prove
Proposition~\ref{p:localization}, which we need for the proof of
Theorem~\ref{t:equivalence}.

Let $A_1$ be an (unbounded) self-adjoint operator in a Hilbert
space ${\mathcal H}_1$ with the domain $\Dom (A_1)$. We assume
that 
$A_1$ is semi-bounded from 
below:
\begin{equation}\label{e:03}
(A_1u,u)\geq 
\lambda_0\|u\|^2,\quad u\in \Dom (A_1)
\end{equation}
with some 
$\lambda_0\leq 0$.

Let $J$ be a self-adjoint bounded operator in 
${\mathcal H}_1$
that maps the domain of $A_1$ into itself, $J: \Dom 
(A_1)\to \Dom
(A_1)$. We assume that $0\leq J\leq \id_{{\mathcal 
H}_1}$.
Introduce a self-adjoint positive bounded operator $J'$ 
in
${\mathcal H}_1$ by the formula $J^2+(J')^2=\id_{{\mathcal 
H}_1}$.
We assume that $J'$ maps the domain of $A_1$ into itself, 
the
operators $[J,[J,A_1]]$ and $[J',[J',A_1]]$ extend to 
bounded
operators in ${\mathcal H}_1$ 
and
\begin{equation}\label{e:05}
\max (\|[J,[J,A_1]]\|,\, 
\|[J',[J',A_1]]\|)\leq \gamma.
\end{equation}

Finally, we assume 
that
\begin{equation}\label{e:04}
(A_1J'u,J'u)\geq \alpha 
\|J'u\|^2,\quad u\in \Dom(A_1)
\end{equation}
for some 
$\alpha>0$.

Denote by $E(\lambda)$ the spectral projection of the 
operator
$A_1$, corresponding to the semi-axis $(-\infty,\lambda]$. 
We have
\begin{equation}\label{e:01}
(A_1E(\lambda)u,E(\lambda)u)\leq 
\lambda\|E(\lambda)u\|^2, \quad
u\in \Dom 
(A_1).
\end{equation}

\begin{prop}\label{p:localization}
If 
$\alpha>\lambda+\gamma$, then we have the following 
estimate
\begin{equation}\label{e:06}
\|JE(\lambda)u\|^2\geq
\frac{\alpha-\lambda-\gamma}{\alpha-\lambda_0}\|E(\lambda)u\|^2,
\quad 
u\in {\mathcal H}_1.
\end{equation}
\end{prop}

\begin{rem}
Note that 
in the case $\lambda<\lambda_0$ the statement is
trivial. In the 
opposite case $\lambda\geq\lambda_0$, since
$\alpha>\lambda+\gamma$ 
and $\gamma\geq 0$, the coefficient,
entering in the right-hand side 
of the formula (\ref{e:06}),
satisfies the 
estimate
\[
0<\frac{\alpha-\lambda-\gamma}{\alpha-\lambda_0}\leq 
1.
\]
\end{rem}

\begin{proof}
We can assume that 
$\lambda\geq\lambda_0$. By the IMS localization
formula (see 
\cite{Sh} and references there), we 
have
\begin{equation}\label{e:07}
A_1=JA_1J+J'A_1J'+\frac{1}{2}[J,[J,A_1]]+\frac{1}{2}[J',[J',A_1]].
\end{equation}
Applying 
this formula to $E(\lambda)u, u\in \Dom (A_1)$, we 
get
\begin{equation}\label{e:08}
\begin{split}
(A_1E(\lambda)u, 
E(\lambda)u)&=(A_1JE(\lambda)u, JE(\lambda)u) +
(A_1J'E(\lambda)u, 
J'E(\lambda)u)\\
&+\frac{1}{2}([J,[J,A_1]]E(\lambda)u,E(\lambda)u)+
\frac{1}{2}([J',[J',A_1]]E(\lambda)u, 
E(\lambda)u).
\end{split}
\end{equation}
Combining \eqref{e:04}, 
\eqref{e:08}, \eqref{e:01}, \eqref{e:03},
\eqref{e:05}, we 
get
\begin{equation*}
\begin{split}
\|J'E(\lambda)u\|^2& \leq 
\frac{1}{\alpha}(A_1J'E(\lambda)u,
J'E(\lambda)u)\\ &= 
\frac{1}{\alpha}((A_1E(\lambda)u,
E(\lambda)u)-(A_1JE(\lambda)u, 
JE(\lambda)u)
-\frac{1}{2}([J,[J,A_1]]E(\lambda)u,E(\lambda)u)\\
&-\frac{1}{2}([J',[J',A_1]]E(\lambda)u, 
E(\lambda)u))\\ 
&\leq
\frac{1}{\alpha}\left((\lambda+\gamma)\|E(\lambda)u\|^2
-\lambda_0\|JE(\lambda)u\|^2 
\right).
\end{split}
\end{equation*}
Hence, we 
have
\begin{equation*}
\|JE(\lambda)u\|^2=\|E(\lambda)u\|^2-\|J'E(\lambda)u\|^2 
\geq
\left(1-\frac{\lambda+\gamma}{\alpha}\right)\|E(\lambda)u\|^2
+\frac{\lambda_0}{\alpha}\|JE(\lambda)u\|^2,
\end{equation*}
that 
immediately implies the required estimate 
\eqref{e:06}.
\end{proof}

\begin{cor}
If $\alpha>\lambda+\gamma$, 
then we have the following 
estimate:
\begin{equation}\label{e:010}
\|J'E(\lambda)u\|^2\leq
\frac{\lambda+\gamma-\lambda_0}{\alpha-\lambda_0}\|E(\lambda)u\|^2,
\quad 
u\in {\mathcal H}_1.
\end{equation}
\end{cor}

\begin{proof}
This 
follows immediately from the equality
$\|Jv\|^2+\|J'v\|^2=\|v\|^2$ 
for any $v\in {\mathcal H}_1$.
\end{proof}

\begin{cor}
If 
$\alpha>\lambda+\gamma$, then we have the following 
estimate
\begin{equation}\label{e:011}
(A_1JE(\lambda)u, 
JE(\lambda)u)\leq
\left(\lambda+\gamma-\lambda_0\frac{\lambda+\gamma-\lambda_0}{\alpha-\lambda_0}
\right)\|E(\lambda)u\|^2,\quad 
u\in \Dom (A_1).
\end{equation}
\end{cor}
\begin{proof}
  From 
\eqref{e:08} and \eqref{e:010}, we 
get
\begin{equation*}
\begin{split}
(A_1JE(\lambda)u, 
JE(\lambda)u)&=(A_1E(\lambda)u,
E(\lambda)u)-(A_1J'E(\lambda)u, 
J'E(\lambda)u)\\
&-\frac{1}{2}([J,[J,A_1]]E(\lambda)u,E(\lambda)u)-\frac{1}{2}([J',[J',A_1]]E(\lambda)u,
E(\lambda)u)\\ 
&\leq
\left((\lambda+\gamma)\|E(\lambda)u\|^2-\lambda_0\|J'E(\lambda)u\|^2
\right)\\ 
&\leq
\left(\lambda+\gamma-\lambda_0\frac{\lambda+\gamma-\lambda_0}{\alpha-\lambda_0}
\right)\|E(\lambda)u\|^2,
\end{split}
\end{equation*}
as 
desired.
\end{proof}

\subsection{Proof of 
Theorem~\ref{t:equivalence}}
In this Section, we will use the 
notation of
Section~\ref{s:abstract-equiv}. We start with the 
following

\begin{prop}\label{p:closedimage}
If 
$\alpha_1>\lambda_1+\gamma_1$ and $$ \lambda_2>\rho\left[
\beta_1 
\left(\lambda_1+\gamma_1+
\frac{(\lambda_1+\gamma_1-\lambda_{01})^2}{\alpha_1-\lambda_1-\gamma_1}\right)+
\varepsilon_1\right], 
$$ then there exists $\varepsilon_0 >0$ 
such
that
\begin{equation}\label{e:closed}
\|E_2(\lambda_2)i_2Jp_1E_1(\lambda_1)u\|_2^2\geq 
\varepsilon_0
\|E_1(\lambda_1)u\|_1^2, \quad u\in 
H.
\end{equation}
\end{prop}

\begin{proof} Applying \eqref{e:16} to 
a function
$p_1E_1(\lambda_1)u, u\in \Dom (A_1)$ and taking into 
account that
$J_1=i_1Jp_1$, we 
get
\begin{equation}\label{e:7}
(A_2i_2Jp_1E_1(\lambda_1)u,i_2Jp_1E_1(\lambda_1)u)_2\leq 
\beta_1
(A_1J_1E_1(\lambda_1)u,J_1E_1(\lambda_1)u)_1+\varepsilon_1
\|J_1E_1(\lambda_1)u\|_1^2.
\end{equation}

Clearly, 
for any $\lambda$ and $l=1,2$ we have the estimate
\begin{equation} 
\label{e:4}
(A_l(\id_{{\mathcal 
H}_l}-E_l(\lambda))u,(\id_{{\mathcal
H}_l}-E_l(\lambda))u)_l \geq 
\lambda\|(\id_{{\mathcal
H}_l}-E_l(\lambda))u\|_l^2, \quad u\in \Dom 
(A_l).
\end{equation}

By \eqref{e:4}, \eqref{e:6} and \eqref{e:rho}, 
it follows 
that
\begin{multline}\label{e:21}
(A_2i_2Jp_1E_1(\lambda_1)u,i_2Jp_1E_1(\lambda_1)u)_2
=(A_2E_2(\lambda_2)i_2Jp_1E_1(\lambda_1)u,E_2(\lambda_2)i_2Jp_1E_1(\lambda_1)u)_2\\
+(A_2(\id_{{\mathcal 
H}_2}-E_2(\lambda_2))i_2Jp_1E_1(\lambda_1)u,
(\id_{{\mathcal 
H}_2}-E_2(\lambda_2))i_2Jp_1E_1(\lambda_1)u)_2\\
\geq 
\lambda_{02}\|E_2(\lambda_2)i_2Jp_1E_1(\lambda_1)u\|_2^2+
\lambda_2\|(\id_{{\mathcal
H}_2}-E_2(\lambda_2))i_2Jp_1E_1(\lambda_1)u\|_2^2\\ 
=
\lambda_2\|i_2Jp_1E_1(\lambda_1)u\|_2^2 
-
(\lambda_2-\lambda_{02})\|E_2(\lambda_2)i_2Jp_1E_1(\lambda_1)u\|_2^2\\
\geq \lambda_2\rho^{-1}\|J_1E_1(\lambda_1)u\|_1^2 
-
(\lambda_2-\lambda_{02})\|E_2(\lambda_2)i_2Jp_1E_1(\lambda_1)u\|_2^2.
\end{multline}

 From the other side, by \eqref{e:7}, \eqref{e:011}, we 
have
\begin{multline}\label{e:22}
(A_2i_2Jp_1E_1(\lambda_1)u,i_2Jp_1E_1(\lambda_1)u)_2\\ 
\leq
\beta_1 
(A_1J_1E_1(\lambda_1)u,J_1E_1(\lambda_1)u)_1+\varepsilon_1
\|J_1E_1(\lambda_1)u\|_1^2\\ 
\leq 
\beta_1
\left(\lambda_1+\gamma_1-\lambda_{01}
\frac{\lambda_1+\gamma_1-\lambda_{01}}{\alpha_1-\lambda_{01}}\right)
\|E_1(\lambda_1)u\|_1^2 
+\varepsilon_1 \|J_1E_1(\lambda_1)u\|_1^2.
\end{multline}

Combining 
\eqref{e:21} and \eqref{e:22}, we 
get
\begin{multline*}
\lambda_2\rho^{-1}\|J_1E_1(\lambda_1)u\|_1^2 
-
(\lambda_2-\lambda_{02})\|E_2(\lambda_2)i_2Jp_1E_1(\lambda_1)u\|_2^2\\
\leq 
\beta_1 
\left(\lambda_1+\gamma_1-\lambda_{01}
\frac{\lambda_1+\gamma_1-\lambda_{01}}{\alpha_1-\lambda_{01}}\right)
\|E_1(\lambda_1)u\|_1^2 
+\varepsilon_1 \|J_1E_1(\lambda_1)u\|_1^2,
\end{multline*}
that 
implies, using 
\eqref{e:06},
\begin{multline*}
\|E_2(\lambda_2)i_2Jp_1E_1(\lambda_1)u\|_2^2 
\\ \geq
\frac{1}{\lambda_2-\lambda_{02}}
\left[(\lambda_2\rho^{-1}-\varepsilon_1)\|J_1E_1(\lambda_1)u\|_1^2
- \beta_1 
\left(\lambda_1+\gamma_1-\lambda_{01}
\frac{\lambda_1+\gamma_1-\lambda_{01}}{\alpha_1-\lambda_{01}}\right)
\|E_1(\lambda_1)u\|_1^2\right] 
\\ \geq
\frac{1}{\lambda_2-\lambda_{02}} 
\left[(\lambda_2\rho^{-1}-\varepsilon_1)
\frac{\alpha_1-\lambda_1-\gamma_1}{\alpha_1-\lambda_{01}}
       - \beta_1 \left(\lambda_1+\gamma_1-\lambda_{01}
\frac{\lambda_1+\gamma_1-\lambda_{01}}{\alpha_1-\lambda_{01}}\right)\right]
\|E_1(\lambda_1)u\|_1^2
\end{multline*}
as desired.
\end{proof}

\begin{rem}
Note that we only used estimate~(\ref{e:16}) (but not (\ref{e:A1A2}))
in the proof of Proposition~\ref{p:closedimage}.
\end{rem}

\begin{proof}[Proof of Theorem~\ref{t:equivalence}]
As above, we will use the notation of
Section~\ref{s:abstract-equiv}. Take arbitrary $\lambda_1\in
(a_1,b_1)$ and $\lambda_2\in (a_2,b_2)$. Consider the bounded
operator $T={\mathcal
V}_2E_2(\lambda_2)i_2Jp_1E_1(\lambda_1){\mathcal V}^{-1}_1$ in
$\chh$. By assumption, $T$ belongs to the von Neumann algebra
$\pi({\mathfrak A})''$.

Since $$ \lambda_2> a_2=\rho\left[ \beta_1 \left(a_1+\gamma_1+
\frac{(a_1+\gamma_1-\lambda_{01})^2}{\alpha_1-a_1-\gamma_1}\right)+
\varepsilon_1\right] $$ and (see Remark~\ref{b1b2})
\begin{align*}
b_1&=\rho\left[ \beta_2 \left(b_2+\gamma_2+
\frac{(b_2+\gamma_2-\lambda_{02})^2}{\alpha_2-b_2-\gamma_2}\right)+
\varepsilon_2\right]\\ & >\rho\left[ \beta_2
\left(\lambda_2+\gamma_2+
\frac{(\lambda_2+\gamma_2-\lambda_{02})^2}{\alpha_2-\lambda_2-\gamma_2}\right)+
\varepsilon_2\right],
\end{align*}
it follows from Proposition~\ref{p:closedimage} that the map $$
E_2(\lambda_2)i_2Jp_1E_1(\lambda_1) =
E_2(\lambda_2)i_2Jp_1E_1(a_1+0) : \Im E_1(\lambda_1)\to \Im
E_2(\lambda_2) $$ is injective and has closed image, and the map
$$ (E_2(\lambda_2)i_2Jp_1E_1(\lambda_1))^* =
E_1(\lambda_1)i_1Jp_2E_2(\lambda_2) =
E_1(b_1-0)i_1Jp_2E_2(\lambda_2): \Im E_2(\lambda_2)\to \Im
E_1(\lambda_1). $$ is injective. Hence, the map
$E_2(\lambda_2)i_2Jp_1E_1(\lambda_1):\Im E_1(\lambda_1)\to \Im
E_2(\lambda_2)$ is bijective.

Let $T=US, U, S\in \pi({\mathfrak A})'',$ be the polar
decomposition of $T$. Since $\Ker T = {\mathcal V}_1(\Im
E_1(\lambda_1))= \Im {\mathcal V}_1E_1(\lambda_1){\mathcal
V}^{-1}_1$ and $\Im T= {\mathcal V}_2(\Im E_2(\lambda_2))=\Im
{\mathcal V}_2E_2(\lambda_2){\mathcal V}^{-1}_2$, $U$ is a partial
isometry that performs the Murray-von Neumann equivalence of the
projections ${\mathcal V}_1E_1(\lambda_1){\mathcal V}^{-1}_1$ and
${\mathcal V}_2E_2(\lambda_2){\mathcal V}^{-1}_2$ in the von
Neumann algebra $\pi({\mathfrak A})''$.

Since the interval $(a_1,b_1)$ does not intersect with the
spectrum of $A_1$, the spectral density function $\tau({\mathcal
V}_1E_1(\lambda_1) {\mathcal V}^{-1}_1)$ is constant for any
$\lambda_1\in (a_1,b_1)$. Using the Murray-von Neumann equivalence
of ${\mathcal V}_1E_1(\lambda_1){\mathcal V}^{-1}_1$ and
${\mathcal V}_2E_2(\lambda_2){\mathcal V}^{-1}_2$ and the tracial
property, we conclude that the spectral density function
$\tau({\mathcal V}_2E_2(\lambda_2){\mathcal V}^{-1}_2)$ is
constant for any $\lambda_2\in (a_2,b_2)$. Since the trace $\tau$
is faithful, the interval $(a_2,b_2)$ does not intersect with the
spectrum of $A_2$, that completes the proof of the first part of
Theorem~\ref{t:equivalence}.

Note that $E_1(\lambda_1) = \chi_{[e^{-t\lambda_1},
\infty)}\left(e^{-t A_1}\right)$. Using the fact that $\lambda_1$
belongs to a gap in the spectrum $A_1$ and
$e^{-tA_1}\in\pi_1({\mathfrak A})$, one can replace
$\chi_{[e^{-t\lambda_1}, \infty)}$ by a continuous function and
obtain that $E_1(\lambda_1)\in \pi_1({\mathfrak A})$ for any
$\lambda_1\in (a_1,b_1)$. Similarly, $E_2(\lambda_2)\in
\pi_2({\mathfrak A})$ for any $\lambda_2\in (a_2,b_2)$. By
assumption, this implies that $T$ belongs to the $C^*$-algebra
$\pi({\mathfrak A})$.

\begin{lemma}\label{l:polar1}
Let ${\mathfrak A}$ be a $C^*$-algebra, $\chh$ a Hilbert space
equipped with a faithful $\ast$-representation of ${\mathfrak A}$,
$\pi: {\mathfrak A}\to \ck(\chh)$. If $P\in \pi({\mathfrak A})$
has closed image and $P=US$ is its polar decomposition, then $U,
S\in \pi({\mathfrak A})$.
\end{lemma}

\begin{proof} We will identify $\mathfrak A$ with $\pi(\mathfrak A)$.
Since $P$ has closed image, the operators $P^*$ and $P^*P$ also
have closed image, hence $0$ is an isolated point in the spectrum
of $P^*P$. Then it follows that the projection on the kernel of
$P$ (or, which is the same, the kernel of $P^*P$) is in the
$C^*$-algebra ${\mathfrak A}$. Clearly $S=\sqrt{P^*P}$ is in
${\mathfrak A}$ and has $0$ its isolated point in the spectrum.
Now we define a bounded operator $S^{(-1)}$ in $\chh$ by the
equalities $S^{(-1)}(Sx)=x, x\in (\Ker S)^{\bot},$ on $\Im S$ and
$S^{(-1)}(x)=0$ on the orthogonal complement $(\Im S)^{\bot}$.
Since $S^{(-1)}$ is given as $f(S)$, where $f$ is a continuous
function on the spectrum of $S$ defined as
$f(\lambda)=\lambda^{-1}$ if $\lambda\not=0$, $f(0)=0$, we see
that $S^{(-1)}$ is in ${\mathfrak A}$. It remains to notice that
$U=PS^{(-1)}$.
\end{proof}

Applying Lemma~\ref{l:polar1} to the operator $T$, we get that the
partial isometry $U\in \pi({\mathfrak A})$ performs the desired
Murray-von Neumann equivalence of the projections ${\mathcal
V}_1E_1(\lambda_1){\mathcal V}^{-1}_1$ and ${\mathcal
V}_2E_2(\lambda_2){\mathcal V}^{-1}_2$ in the $C^*$-algebra
$\pi({\mathfrak A})$.
\end{proof}

\begin{cor}\label{c:K}
Under assumptions of Theorem~\ref{t:equivalence}, for any
$\lambda_1\in (a_1,b_1)$ and $\lambda_2\in (a_2,b_2)$, the
spectral projections ${\mathcal V}_1E_1(\lambda_1){\mathcal
V}^{-1}_1$ and ${\mathcal V}_2E_2(\lambda_2){\mathcal V}^{-1}_2$
define the same element in K-theory of ${\mathfrak A}$:
\[
[{\mathcal V}_1E_1(\lambda_1){\mathcal V}^{-1}_1]=[{\mathcal
V}_2E_2(\lambda_2){\mathcal V}^{-1}_2]\in K_0(\mathfrak A).
\]
\end{cor}

\subsection{The model operator and equivalence of spectral
projections in the $C^*$-algebra} \label{s:model1} As above, let
$M$ be a compact connected Riemannian manifold, $\Gamma$ its
fundamental group and $\widetilde M$ its universal cover. Let $\E$
be a Hermitian vector bundle on $M$ and $\tE$ the lift of $\E$ to
the universal cover $\widetilde M$. Let $\omega$ be a closed
real-valued 2-form on $M$ such that ${\bf B} =p^* \omega$ is
exact, ${\bf B}=d{\bf A}$, where ${\bf A}$ is a real-valued 1-form
on $\widetilde M$.

As above, consider the elliptic self-adjoint differential operator
in $L^2(\M,\tE)$ given by $$ H(\mu) = \mu \nabla^*\nabla + B  +
\mu^{-1} V, $$ where $\nabla$ is a Hermitian connection on the
vector bundle $\tE$ over $\widetilde{M}$ of the form $\nabla =
\widetilde\nabla^\E\otimes \id + \id\otimes \nabla_{\bf A}$, where
$\widetilde \nabla^\E$ is a $\Gamma$-invariant Hermitian
connection on $\tE$, $\nabla_{\bf A}=d+\,i{\bf A}$ is a Hermitian
connection on the trivial vector bundle, $B$, $V$ are
$\Gamma$-invariant self-adjoint endomorphisms of the bundle $\tE$,
and where $V$ satisfies in addition the following {\em Morse type
condition:} $V(x) \ge 0$ for all $x\in \M$. Also if the matrix
$V(x_0)$ is degenerate for some $x_0$ in $\M$, then $V(x_0)=0$ and
there is a positive constant $c$ such that $V(x)\ge c|x-x_0|^2 I$
for all $x$ in a neighborhood of $x_0$. We will also assume that
$V$ has at least one zero point.

Choose a fundamental domain $\F\subset\widetilde M$ so that there
is no zeros of $V$ on the boundary of $\F$. This is equivalent to
saying that  the translations $\{\gamma\F,\;\gamma\in\Gamma\}$
cover  the set $V^{-1}(0)$ (the set of all zeros of $V$). Let
$V^{-1}(0)\cap\F=\{\bar x_j|\,j=1,\dots,N\}$ be the set of all
zeros of $V$ in $\F$; $\bar x_i\ne\bar x_j$ if $i\ne j$.

Let $K$ denote the {\em model operator} of $H$ (cf. \cite{Sh}),
which is obtained as a direct sum of quadratic parts of $H$ in all
points $\bar x_1,\dots,\bar x_N$. More precisely, $K$ is an
operator in $L^2({\mathbb R}^n,\C^k)^N$ given by $$ K =
\oplus_{1\le j\le N} K_j, $$ where $K_j$ is an unbounded
self-adjoint operator in $L^2({\mathbb R}^n,\C^k)$ which
corresponds to the zero $\bar x_j$. It is a quantum harmonic
oscillator and has a discrete spectrum. We assume that we have
fixed local coordinates on $\M$ and trivialization of the bundle
$\tE$ in a small neighborhood $B(\bar x_j, r)$ of $\bar x_j$  for
every $j=1,\dots,N$. We assume that $\bar x_j$ becomes zero in
these local coordinates. Then $K_j$ has the form $$ K_{j}=
H_{j}^{(2)}+\bar{B}_j+ V^{(2)}_{j}, $$ where all the components
are obtained from $H$ as follows. The second order term
$H_{j}^{(2)}$ is a homogeneous second order  differential operator
with constant coefficients (without lower order terms) obtained by
isolating the second order terms in the operator $H$ and freezing
the coefficients of this operator at $\bar x_j$. (Note that
$H_j^{(2)}$ does not depend on $\mathbf A$.) The zeroth order term
$V^{(2)}_{j}$ is  obtained by taking the quadratic part of $V$ in
the chosen coordinates near $\bar x_j$.

More explicitly, $$ H_j^{(2)}=-\sum_{i,k=1}^n g^{ik}(\bar
x_j)\frac{\partial^2}{\partial x_i\partial x_k}, \qquad V_j^{(2)}=
\frac{1}{2}\sum_{i,k=1}^n\frac{\partial^2 V}{\partial x_i\partial
x_k}(\bar x_j)x_ix_k, $$ where $(g^{ik})$ is the inverse matrix to
the matrix of the Riemannian tensor $(g_{ik})$.

Finally, $\bar{B}_j=B(\bx_j),\ j=1,\dots, N$, so $\bar{B}_j$ is an
endomorphism of the fiber of the bundle $\tE$ over the point
$\bx_j$.

We will also need the operator $$ K(\mu) = \oplus_{1\le j\le N}
K_j(\mu), $$ where $$ K_j(\mu)=\mu H_{j}^{(2)}+\bar{B}_j+\mu^{-1}
V^{(2)}_{j}, \qquad \mu>0. $$ It is easy to see that $K(\mu)$ has
the same spectrum as $K=K(1)$.

We will say that $H$ is {\it flat} near $\bar x_j$ if
$H(\mu)=K_j(\mu)$ for all $\mu$ near $\bar x_j$. (In particular,
in this case we should have $A=0$ near $\bar x_j$.)

We are going to apply Theorem~\ref{t:equivalence} in the following
particular setting. Take the $C^*$ algebra ${\mathfrak A}$ to be
$C^*_r(\Gamma, \bar\sigma)\otimes {\mathcal K}$. Let $\mathcal H$
be the Hilbert space $\ell^2(\Gamma)\otimes \ell^2(\N)$. Put
${\mathcal H}_1=\ell^2(\Gamma)\otimes L^2({\mathbb R}^n,\C^k)^N$
($\ell^2(\Gamma)\otimes\mathfrak{H}_K$ in the above notation) and
${\mathcal H}_2=L^2(\M,\tE)$ (denoted by $\mathfrak{H}$ above).
Choose an arbitrary unitary isomorphism $V_1:L^2({\mathbb
R}^n,\C^k)^N\to \ell^2(\N)$ and define an unitary operator
${\mathcal V}_1 : {\mathcal H}_1\to {\mathcal H}$ as ${\mathcal
V}_1=\id\otimes V_1$. Similarly, choose an arbitrary unitary
isomorphism $V_2:L^2(\F,\tE|_\F)\to \ell^2(\N)$ and define an
unitary operator ${\mathcal V}_2 : {\mathcal H}_2\to {\mathcal H}$
as ${\mathcal V}_2=(\id\otimes V_2)\circ {\bf U}$, where ${\bf U}$
is the $(\Gamma, \sigma)$-equivariant isometry \eqref{isom}.

Let $\pi$ be the representation of the algebra ${\mathfrak A}$ in
$\mathcal H$ given by the tensor product of the representation of
$C^*_r(\Gamma, \bar\sigma)$ on $\ell^2(\Gamma)$ by left twisted
convolutions and the natural representation of ${\mathcal K}$ in
$\ell^2(\N)$. So we have $\pi(C^*_r(\Gamma, \bar\sigma)\otimes
{\mathcal K})\subset {\mathcal
A}^L_{\ell^2(\N)}(\Gamma,\bar\sigma)$ and $\pi(C^*_r(\Gamma,
\bar\sigma)\otimes {\mathcal K})''= {\mathcal
A}^L_{\ell^2(\N)}(\Gamma,\bar\sigma)$. Using the unitary
isomorphisms ${\mathcal V}_1$ and ${\mathcal V}_2$, we get
representations $\pi_1$ and $\pi_2$ of ${\mathfrak A}$ in
${\mathcal H}_1$ and ${\mathcal H}_2$ accordingly,
$\pi_l(a)={\mathcal V}^{-1}_l\circ\pi(a)\circ {\mathcal V}_l,
l=1,2, a\in {\mathfrak A}$.

Consider self-adjoint, semi-bounded from below operators $A_1=\id
\otimes K(\mu)$ in ${\mathcal H}_1$ and $A_2=H(\mu)$ in ${\mathcal
H}_2$. It is clear that $e^{-tA_1}=\id \otimes e^{-tK(\mu)}\in
\pi_1({\mathfrak A})\cong C^*_r(\Gamma, \bar\sigma)\otimes
{\mathcal K}(\mathfrak{H}_K)$ for any $t>0$. As shown in
Lemma~\ref{F} below, for any $t>0$, the operator $e^{-tA_2}$
belongs to $\pi_2({\mathfrak A})$. Remark that, in notation of
Theorem~\ref{main1}, $E_1(\lambda)=\id\otimes E^0(\lambda)$ and
$E_2(\lambda)=E(\lambda)$.

Let ${\mathcal H}_0=\ell^2(\Gamma)\otimes \left(\oplus_{j=1}^N
L^2(B(\bar x_j, r),\tE|_{B(\bar x_j,r)})\right)$.  An inclusion $
i_1: {\mathcal H}_0 \to {\mathcal H}_1$ is defined as $i_1=\id
\otimes j_1$, where $j_1$ is the inclusion of
$\left(\oplus_{j=1}^NL^2(B(\bar x_j, r),\tE|_{B(\bar
x_j,r)})\right)$ in $L^2({\mathbb R}^n,\C^k)^N$ given by the
chosen local coordinates and trivializations of the vector bundle
$\tE$. An inclusion $i_2: {\mathcal H}_0 \to {\mathcal H}_2$ is
defined as $i_2={\bf U}^{*}\circ (\id\otimes j_2)$, where $j_2$ is
the natural inclusion of $\left(\oplus_{j=1}^NL^2(B(\bar x_j,
r),\tE|_{B(\bar x_j,r)})\right)$ in $L^2(\F, \tE|_\F)$.

The operator $p_1:{\mathcal H}_1\to {\mathcal H}_0$ is defined as
$p_1=\id\otimes r_1$, where $r_1$ is the restriction operator
$L^2({\mathbb R}^n,\C^k)^N\to \oplus_{j=1}^NL^2(B(\bar x_j,
r),\tE|_{B(\bar x_j,r)})$. The operator $p_2:{\mathcal H}_1\to
{\mathcal H}_0$ is defined as $p_2=(\id\otimes r_2)\circ {\bf U}$,
where $r_2: L^2(\F,\tE|_\F)\to \oplus_{j=1}^NL^2(B(\bar x_j,
r),\tE|_{B(\bar x_j,r)})$ is the restriction operator.

Fix a function $\phi\in C_0^\infty(\R^n)$ such that $0\leq\phi\leq
1$, $\phi(x)=1$ if $|x|\leq 1$, $\phi(x)=0$ if $|x|\geq 2$, and
$\phi'=(1-\phi^2)^{1/2}\in C^\infty(\R^n)$. Fix a number $\ka,\
0<\ka<1/2,$ which we shall choose later. For any $\mu>0$ define
$\phi^{(\mu)}(x)=\phi(\mu^{-\ka}x)$. For any $\mu>0$ small enough,
let $\phi_j=\phi^{(\mu)}\in C^\infty_c(B(\bar x_j,r))$ in the
fixed coordinates near $\bar x_j$. Denote also $\phi_{j,\gamma} =
(\gamma^{-1})^*\phi_j$. (This function is supported near
$\gamma\bar x_j$.) We will always take $\mu\in (0,\mu_0)$ where
$\mu_0$ is sufficiently small, so in particular the supports of
all functions $\phi_{j,\gamma}$ are disjoint.

Let $\Phi=\oplus_{j=1}^N \phi_j\in \oplus_{j=1}^NC^\infty_c(B(\bar
x_j,r))\subset C^\infty(\F)$. Consider a $(\Gamma,
\sigma)$-equivariant, self-adjoint, bounded operator $J$ in
${\mathcal H}_0$ defined as $J=\id\otimes \Phi$, where $\Phi$
denotes the multiplication operator by the function $\Phi$ in the
space $\oplus_{j=1}^NL^2(B(\bar x_j,r),\tE|_{B(\bar x_j,r)})$.

It is clear that ${\mathcal V}_2i_2Jp_1{\mathcal
V}_1^{-1}=\id\otimes V_2j_2J_0r_1V_1^{-1}$, where $j_2J_0r_1$ is
the multiplication operator by the function $\Phi$, considered as
an operator from $L^2({\mathbb R}^n,\C^k)^N$ to $L^2(\F,\tE|_\F)$.
Hence, one can easily see that the operator ${\mathcal
V}_2i_2Jp_1{\mathcal V}^{-1}_1$ belongs to the von Neumann algebra
$\pi({\mathfrak A})'' = {\mathcal A}^L_{\ell^2(\N)}
(\Gamma,\bar\sigma)\cong {\mathcal A}^L (\Gamma,\bar\sigma)\otimes
{\mathcal B}(\ell^2(\N))$, $(i_2Jp_1)^*=i_1Jp_2$ and, for any
$a\in {\mathfrak A}$, the operator $\pi(a){\mathcal
V}_2(i_2Jp_1){\mathcal V}^{-1}_1$ belongs to $\pi({\mathfrak A})$.

We will use local coordinates near $\bar{x}_j$ such that the
Riemannian volume element at the point $\bar{x}_j$ coincides with
the Euclidean volume element given by the chosen local
coordinates. Similarly we will fix a trivialization of the bundle
$\tE$ near $\bar{x}_j$ such that the Hermitian metric becomes
trivial in this trivialization. Then the estimate \eqref{e:rho}
holds with $\rho=1+O(\mu^\kappa)$.

Denote by the same letters $\phi$ and $\phi'$ the multiplication
operators in $L^2({\mathbb R}^n,\C^k)$ by the functions $\phi$ and
$\phi'$ accordingly. Let $\Phi_1$ and $\Phi'_1$ be the bounded
operators in $L^2({\mathbb R}^n,\C^k)^N\cong L^2({\mathbb
R}^n,\C^k)\otimes \C^N$ given by $\Phi_1=\phi\otimes \id_{\C^N}$
and $\Phi'_1=\phi'\otimes \id_{\C^N}$. Then we have
$J_1=\id\otimes \Phi_1 $ and $J'_1=\id\otimes \Phi'_1$ in
$\ell^2(\Gamma)\otimes L^2({\mathbb R}^n,\C^k)^N$. Let
$\Phi_\gamma=\oplus_{j=1}^N \phi_{j,\gamma}\in
\oplus_{j=1}^NC^\infty_c(B(\gamma\bar x_j,r))\subset
C^\infty(\widetilde M)$ and
$\Phi_2=\sum_{\gamma\in\Gamma}\Phi_\gamma \in
C^{\infty}(\widetilde M)$. Define a function $\Phi'_2\in
C^{\infty}(\widetilde M), \Phi'_2\geq 0$ by the equation
$(\Phi_2)^2+(\Phi'_2)^2=1\ \text{in}\ C^{\infty}(\widetilde M)$.
The operators $J_2$ and $J'_2$ are given by the multiplication
operators by the functions $\Phi_2$ and $\Phi'_2$ in
$L^2(\widetilde{M},\tE)$ respectively.

Let $a_{1,j}^{(2)}, j=1,2,\ldots,N,$ be the principal symbol of
$K_j$, which is a function on $T^*{\mathbb R}^n$:
\[
a_{1,j}^{(2)}(x,\xi)=\sum_{i,k=1}^n g^{ik}(\bar x_j)\xi_i\xi_k,
\quad (x,\xi)\in T^*{\mathbb R}^n.
\]
Then the operators $[J_1,[J_1,A_1]], [J'_1,[J'_1,A_1]]$ are given
by $ -\mu \id \otimes \left( \oplus_{1\leq j\leq N}
a_{1,j}^{(2)}(x,d\phi(x))\right)$ and $ -\mu \id \otimes
\left(\oplus_{1\leq j\leq N} a_{1,j}^{(2)}(x,d\phi'(x))\right)$ in
$\ell^2(\Gamma)\otimes L^2({\mathbb R}^n,\C^k)^N$ accordingly.
Similarly, let $a_2^{(2)}$ be the principal symbol of $H(1)$,
which is a function on $T^*\widetilde M$:
\[
a_{2}^{(2)}(x,\xi)=\sum_{i,k=1}^n g^{ik}(x)\xi_i\xi_k, \quad
(x,\xi)\in T^*\widetilde M.
\]
Then $[J_2,[J_2,A_2]], [J'_2,[J'_2,A_2]]$ are the multiplication
operators by functions $-\mu a_2^{(2)}(x,d\Phi_2(x))$ and $-\mu
a_2^{(2)}(x,d\Phi'_2(x))$ in $L^2(\widetilde M, \tE)$. Therefore,
\begin{align*}
\gamma_1&=\mu \max_{j=1,2,\ldots,N} \max \left( \sup_{ x\in
{\mathbb R}^n} (a^{(2)}_{1,j}(x,d\phi(x))), \sup_{x\in {\mathbb
R}^n} (a^{(2)}_{1,j}(x,d\phi'(x)))\right) =O(\mu^{1-2\kappa}),\\
\gamma_2&=\mu \max \left(\sup_{x\in \M} (a^{(2)}_2(x,d\Phi_2(x))),
\sup_{x\in \M} (a^{(2)}_2(x,d\Phi'_2(x)))\right)
=O(\mu^{1-2\kappa}).
\end{align*}
Since there exists $c_0>0$ such that $V_j^{(2)}\ge
c_0\mu^{2\kappa}, j=1,2,\ldots,N,$ on $\supp \phi'$ and $V\ge
c_0\mu^{2\kappa}$ on $\supp \Phi'_2$, the estimates \eqref{e:14}
hold with $\alpha_l=c\mu^{-1+2\kappa}$ (see \cite[Lemma 3.3]{Sh}
for more details).

The constants $\lambda_{0l}, l=1,2,$ can be chosen to be
independent of $\mu$. One can take
\[
\lambda_{01}=\lambda_0(K(1))_-,
\]
where $\lambda_0(K(1))$ is the bottom of the spectrum of the
operator $K(1)$ in $L^2({\mathbb R}^n,\C^k)^N$ ($a_-=\min(a,0)$)
and
\[
\lambda_{02}=\inf_{x\in \widetilde M}\lambda_{0}(B(x))_-,
\]
where $\lambda_{0}(B(x)), x\in \M,$ is the lowest eigenvalue of
the linear map $B(x):\tE_x\to\tE_x$ given by the action of the
endomorphism $B$ in fibres of $\tE$.

Finally, the estimates \eqref{e:16} and \eqref{e:A1A2} hold with
$\beta_l=1+O(\mu^\kappa)$ and $\varepsilon_l=O(\mu^{3\kappa-1})$.
This is an immediate consequence of the next lemma, which is an
easy extension of Lemma 3.4 in \cite{Sh}. We will state this lemma
in a slightly more general situation than we need in this paper.

Let $B(0,r)$ denote the open ball in $\R^n$ with radius $r$
centered at the origin. Consider volume elements $\omega_1$ and
$\omega_2$ of the form $\omega_l=\sqrt{g_l(x)}\,dx, l=1,2,$ where
$g_l\in C^\infty(B(0,r))$ and $g_l>0$. Let $(\cdot,\cdot)_l,
l=1,2,$ denote the inner products in $C^\infty_c(B(0,r),\C^k)$
given by the volume elements $\omega_l$ and the standard Hermitian
structure in $\C^k$ and $L^2(B(0,r),\C^k,\omega_l)$ the Hilbert
space of square integrable $\C^k$-valued functions on $B(0,r)$
equipped with the inner product $(\cdot,\cdot)_l$.

Consider formally self-adjoint differential operators $T_l,
l=1,2,$ in $L^2(B(0,r),\C^k,\omega_l)$, depending on a small
parameter $\mu>0$, of the form
\[T_l=-\mu D_l+B_l+\mu^{-1}V_l,\quad l=1,2.\] Here $D_l$ is a second
order formally self-adjoint uniformly elliptic differential
operator with a negative principal symbol, so that $-D_l$ is
semi-bounded from below on $C^\infty_c(B(0,r),\C^k)$; $B_l$ and
$V_l$ are zero order formally self-adjoint operators, i.e. the
multiplication operators by Hermitian $k\times k$ matrix functions
$B_l$ and $V_l$ respectively.

Being formally self-adjoint the operator $D_l$ should have the
form $$D_l=A_l^{(2)}+A_l^{(1)}+A_l^{(0)},$$ where $A_l^{(s)}$ is
an operator of order $s$, $s=0,1,2$,
\begin{align*}
A_l^{(2)}&=\sum_{1\leq r,s\leq n}
{\frac{1}{\sqrt{g_l}}}\frac{\pa}{\pa x^r}
\sqrt{g_l}A_{l,rs}(x){\frac{\pa}{\pa x^s}},\ \
A_{l,rs}^\ast=A_{l,sr};\\ A_l^{(1)}&=\sum_{1\leq r\leq
n}A_{l,r}(x){\frac{\pa}{\pa x^r}},
\end{align*}
$A_{l,rs}$ and $A_{l,r}$ are $k\times k$ smooth matrix functions
on $B(0,r)$; $A_l^{(0)}$ is just a multiplication by a smooth
matrix function $A_l^{(0)}(x)$.

The principal symbol of $-D_l$ is the matrix function on
$B(0,r)\times\R^n$ $$a_l^{(2)}(x,\xi)=\sum_{1\leq r,s\leq
n}\xi_r\xi_s A_{l,rs}(x)\;.$$  For the self-adjoint operator
$-D_l$ its uniform ellipticity and semi-boundedness from below
mean that the matrix $a_l^{(2)}(x,\xi)$ is positive definite for
all $(x,\xi)\in B(0,r)\times\R^n$  and
\begin{equation}\label{e:symbol}
a_l^{(2)}(x,\xi)\geq C_l|\xi |^2, \quad (x,\xi)\in
B(0,r)\times\R^n,
\end{equation}
with some constants $C_l>0$.

Let us assume
\begin{gather*}
A_{1,rs}(0)=A_{2,rs}(0),\quad B_1(0)=B_2(0),\quad
V_1(0)=V_2(0)=0,\\ \frac{\pa V_1}{\pa x^r}(0)=\frac{\pa V_2}{\pa
x^r}(0)=0, \quad r=1,2,\ldots,n,\\ \frac{\pa^2 V_1}{\pa x^r\pa
x^s}(0)=\frac{\pa^2 V_2}{\pa x^r \pa x^s}(0), \quad
r,s=1,2,\ldots,n.
\end{gather*}
Assume also that $g_l(0)=1, l=1,2$, i.e. the volume elements
$\omega_l=\sqrt{g_l(x)}\,dx, l=1,2,$ at the origin coincide with
the Euclidean volume element in $\R^n$.

Finally, as above, let $\phi\in C_c^\infty(\R^n)$ satisfy
$0\leq\phi\leq 1$, $\phi(x)=1$ if $|x|\leq 1$, $\phi(x)=0$ if
$|x|\geq 2$, and, for any $\mu>0$ small enough, define
$\phi^{(\mu)}(x)=\phi(\mu^{-\ka}x), x\in B(0,r)$. Denote by the
same letter $\phi^{(\mu)}$ the multiplication operator in
$C^\infty_c(B(0,r),\C^k)$ by the function $\phi^{(\mu)}$.

\begin{lemma}
Let $1/3<\kappa<1/2$. There exist $C>0$ and $\mu_0>0$ such that
for any $\mu\in (0,\mu_0)$,
\begin{equation}\label{e:T2T1}
(T_2\phi^{(\mu)}u,\phi^{(\mu)}u)_2\leq
(1+C\mu^\ka)(T_1\phi^{(\mu)}u,\phi^{(\mu)}u)_1+C\mu^{3\ka-1}(\phi^{(\mu)}u,\phi^{(\mu)}u)_1
\end{equation}
for any $u\in C^\infty_c(B(0,r),\C^k)$.
\end{lemma}

\begin{proof}
We want to estimate from above the quadratic form
$(T_2\phi^{(\mu)}u,\phi^{(\mu)}u)_2$ in terms of the form
$(T_1\phi^{(\mu)}u,\phi^{(\mu)}u)_1$ for any $u\in
C^\infty_c(B(0,r),\C^k)$. We start with the term
\begin{equation}\label{e:A2}
(-\mu A_2^{(2)}\phi^{(\mu)}u,\phi^{(\mu)}u)_2=\mu\int_{|x|\leq
2\mu^\ka} \sum_{1\leq r,s\leq n} \left(A_{2,rs}(x)\frac{\pa
(\phi^{(\mu)} u)}{\pa x^s},\frac{\pa \overline{(\phi^{(\mu)}
u)}}{\pa x^r}\right) \sqrt{g_2}\,dx.
\end{equation}
Denote by $D(\phi^{(\mu)} u)$ the $k\times n$ matrix of all first
partial derivatives of all components of the vector $\phi^{(\mu)}
u$. Let $\Vert D(\phi^{(\mu)} u)\Vert_0$ be the $L^2$ norm of
$D(\phi^{(\mu)} u)$ considered as a vector function:
\[
\Vert D(\phi^{(\mu)} u)\Vert_0=\left(\int_{|x|\leq 2\mu^\ka}
\sum_{1\leq r\leq n}  \left|\frac{\pa (\phi^{(\mu)} u)}{\pa
x^r}\right|^2\,dx\right)^{1/2}.
\]
Similarly denote by $\Vert \phi^{(\mu)} u\Vert_0$ the $L^2$-norm
of $\phi^{(\mu)} u$ with respect to the Euclidean volume form:
\[
\Vert \phi^{(\mu)} u\Vert_0=\left(\int_{|x|\leq 2\mu^\ka}
\left|\phi^{(\mu)} u\right|^2 \,dx\right)^{1/2}.
\]
Since $g_l(0)=1, l=1,2$, we have
\[
(1-C\mu^\ka) \Vert \phi^{(\mu)} u\Vert_0^2\leq \Vert \phi^{(\mu)}
u\Vert_l^2\leq  (1+C\mu^\ka) \Vert \phi^{(\mu)} u\Vert_0^2, \quad
l=1,2,
\]
with some constant $C>0$. From (\ref{e:symbol}) it follows that
\begin{equation}\label{e:rough}
C_1 \Vert D(\phi^{(\mu)} u)\Vert_0^2 \leq (-
A_l^{(2)}\phi^{(\mu)}u,\phi^{(\mu)}u)_0 \leq C_2 \Vert
D(\phi^{(\mu)} u)\Vert_0^{2},\quad l=1,2,
\end{equation}
with some constants $C_1, C_2>0$. Replacing in (\ref{e:A2})
$\sqrt{g_2}$ by $\sqrt{g_1}$ and $A_{2,rs}(x)$ by $A_{1,rs}(x)$,
we add terms of similar form but with additional factor
$O(\mu^\ka)$. Taking into account (\ref{e:rough}), we get that
there exist $C>0$ and $\mu_0>0$ such that for any $\mu\in
(0,\mu_0)$ $$(-\mu A_2^{(2)}\phi^{(\mu)}u,\phi^{(\mu)}u)_2\leq
(1+C\mu^\ka)(-\mu A_1^{(2)}\phi^{(\mu)}u,\phi^{(\mu)}u)_1.$$

Now estimate the term  $(-\mu
A_2^{(1)}\phi^{(\mu)}u,\phi^{(\mu)}u)_2$. Then we obviously have
for every $\eps>0$ $$ |(-\mu
A_2^{(1)}\phi^{(\mu)}u,\phi^{(\mu)}u)_2|\leq C\mu\Vert
D(\phi^{(\mu)}u)\Vert_0\Vert \phi^{(\mu)}u\Vert_0\leq
C\mu\eps\Vert D(\phi^{(\mu)}u)\Vert_0^2+C\mu\eps^{-1}\Vert
\phi^{(\mu)}u\Vert_0^2. $$ Taking $\eps=\mu^\ka$ and using
(\ref{e:rough})  we obtain $$|(-\mu
A_2^{(1)}\phi^{(\mu)}u,\phi^{(\mu)}u)_2|\leq C\mu^\ka(-\mu
A_1^{(2)}\phi^{(\mu)}u,\phi^{(\mu)}u)_1+
C\mu^{1-\ka}(\phi^{(\mu)}u,\phi^{(\mu)}u)_1,\ \ \mu\in
(0,\mu_0).$$

Also obviously $$|(-\mu
A_2^{(0)}\phi^{(\mu)}u,\phi^{(\mu)}u)_2|\leq C\mu(\phi^{(\mu)}u,
\phi^{(\mu)}u)_1.$$

Therefore we obtain for small $\mu$ $$(-\mu D_2 u,u)\leq
(1+C\mu^\ka)(-\mu
A_1^{(2)}\phi^{(\mu)}u,\phi^{(\mu)}u)_1+C\mu^{1-\ka}(\phi^{(\mu)}u,
\phi^{(\mu)}u)_1.$$

Replacing $B_2(x)$ by $B_1(x)$ and $\sqrt{g_2}$ by $\sqrt{g_1}$ in
the quadratic form $(B_2\phi^{(\mu)}u,\phi^{(\mu)}u)_2$
contributes a term which can be estimated by
$C\mu^\ka(\phi^{(\mu)}u,\phi^{(\mu)}u)_1$:
\[
|(B_2\phi^{(\mu)}u,\phi^{(\mu)}u)_2-(B_1\phi^{(\mu)}u,\phi^{(\mu)}u)_1|\leq
C\mu^\ka(\phi^{(\mu)}u,\phi^{(\mu)}u)_1.
\]

Finally, we have
\[|\mu^{-1}(V_2\phi^{(\mu)}u,\phi^{(\mu)}u)_2-\mu^{-1}(V_1\phi^{(\mu)}u,\phi^{(\mu)}u)_1|\leq
C\mu^{3\ka-1}(\phi^{(\mu)}u,\phi^{(\mu)}u)_1.\] Gathering together
all these estimates, we obtain (\ref{e:T2T1}).
\end{proof}

Now we complete the proofs of Theorem~\ref{main0}, the first part
of Theorem~\ref{main1} and Corollary~\ref{main3}. Assume that the
interval $[a,b]$ does not intersect with the spectrum of $A_1$.
Then there exists an open interval $(a_1,b_1)$ that contains
$[a,b]$ and does not intersect with the spectrum of $A_1$. Using
the formulas
\begin{gather*}
\rho=1+O(\mu^\kappa),\quad \alpha_l=O(\mu^{-1+2\kappa}),\quad
\beta_l=1+O(\mu^\kappa),\\ \varepsilon_l=O(\mu^{3\kappa-1}),\quad
\lambda_{0l}=O(1),\quad \gamma_l=O(\mu^{1-2\kappa}),
\end{gather*}
one can see that, for $a_2$ and $b_2$ given by \eqref{e:a2} and
\eqref{e:b2}, we have
\begin{equation}\label{e:estimate}
a_2=a_1+O(\mu^{s}),\quad b_2=b_1+O(\mu^{s}), \quad \mu\to 0,
\end{equation}
where $s=\min\{3\kappa-1,1-2\kappa\}$. The best possible value of
$s$ which is \[s=\max_\kappa\min\{3\kappa-1,1-2\kappa\}=1/5\] is
attained when $\kappa=2/5$.

Hence, if $\mu>0$ is small enough, we have
$\alpha_1>a_1+\gamma_1$, $\alpha_2>b_2+\gamma_2$, $b_2>a_2$ and
the interval $(a_2,b_2)$ contains $[a,b]$. By
Theorem~\ref{t:equivalence}, we conclude that the interval
$(a_2,b_2)$ does not intersect with the spectrum of $A_2$, that
completes the proof of Theorem~\ref{main0}. Moreover, we have
that, for any $\lambda_1\in (a_1,b_1)$ and $\lambda_2\in
(a_2,b_2)$, the spectral projections ${\mathcal
V}_1E_1(\lambda_1){\mathcal V}^{-1}_1$ and ${\mathcal
V}_2E_2(\lambda_2){\mathcal V}^{-1}_2$ are equivalent in
${\mathfrak A}$. Putting $U={\mathcal V}^{-1}_1{\mathcal V}_2$, we
get the desired Murray - von Neumann equivalence of
$E_1(\lambda_1)=\id\otimes E^0(\lambda)$ and ${\mathcal
V}^{-1}_1{\mathcal V}_2E_2(\lambda_2){\mathcal V}^{-1}_2{\mathcal
V}_1=U E(\lambda) U^{-1}$ in $\pi_1(\mathfrak A)= C^*_r(\Gamma,
\bar\sigma)\otimes \mathcal K(L^2({\mathbb R}^n, \C^k)^N)$. This
immediately implies (\ref{e:1}) (see also Corollary~\ref{c:K}).

Fix an arbitrary isomorphism ${\mathfrak H}_K\cong \ell^2(\N)$. It
induces an isomorphism ${\mathcal K}({\mathfrak H}_K)\cong
{\mathcal K}$ and allows us to write any element $x\in {\mathcal
K}({\mathfrak H}_K)$ as a matrix $(x_{ij}, i,j\in\N)$. Recall that
the Morita equivalence of $K$-theory $K_0(C^*_r(\Gamma,
\bar\sigma)) \cong K_0(C^*_r(\Gamma, \bar\sigma)\otimes \mathcal
K({\mathfrak H}_K)) $  is induced by the standard algebra
homomorphism $\mathcal A \to \mathcal A \otimes \mathcal K$ which
maps $a\in \mathcal A$ to a matrix with the left-upper corner
matrix element $a$, the rest matrix elements being $0$. The
projection $E^0(\lambda)$ belongs to some matrix algebra
$M_n(\C)\subset \mathcal K$, and, under the above isomorphism, the
element $[\id\otimes E^0(\lambda)]\in K_0(C^*_r(\Gamma,
\bar\sigma)\otimes \mathcal K({\mathfrak H}_K))$ corresponds to
the element $[E^0(\lambda)]\in K_0(C^*_r(\Gamma, \bar\sigma))$
given by the matrix $1\otimes E^0(\lambda)\in C^*_r(\Gamma,
\bar\sigma)\otimes M_n(\C)= M_n(C^*_r(\Gamma, \bar\sigma))$. It is
easy to see that the element $[E^0(\lambda)]$ belongs to the image
of $K_0(\C)$ in $K_0(C^*_r(\Gamma, \bar\sigma))$, and the
corresponding class in $\tilde{K}_0(C^*_r(\Gamma, \bar\sigma))$
vanishes, that proves (\ref{e:11}) and completes the proof of the
first part of Theorem~\ref{main1}.

Now Corollary~\ref{main3} follows immediately from the equality
(\ref{e:1}) (see also Lemma~\ref{l:vanishing} below).

\begin{rem}
If $H$ is flat near all points $\bar x_j$, then the estimate
(\ref{e:estimate}) can be improved as follows: for any
$\varepsilon>0$
\[
a_2=a_1+O(\mu^{1-\varepsilon}),\quad
b_2=b_1+O(\mu^{1-\varepsilon}), \quad \mu\to 0.
\]
Indeed, in this case we have (see \cite[Lemma 3.4]{Sh})
\[
\phi_{j}H(\mu)\phi_{j}=\phi_{j}K_j(\mu)\phi_{j},\quad \mu > 0,
\]
therefore,
\[
\beta_l=1,\quad \varepsilon_l=0, \quad l=1,2.
\]
Taking this into account, we easily get
$a_2=a_1+O(\mu^{1-2\kappa}), b_2=b_1+O(\mu^{1-2\kappa})$ as $\mu\to 0$
for any $0<\kappa<1/2$ as desired.
\end{rem}

\section{Review of smooth algebras and higher traces}\label{smooth}

\subsection{Definition and properties of the $*$-algebra $\B(\Gamma,
\sigma)$}

We begin by recalling some generalities on smooth subalgebras of
$C^*$-algebras. Let ${\mathfrak A}$ be a $C^*$-algebra and
$\widetilde {\mathfrak A}$ be obtained by adjoining a unit to
${\mathfrak A}$. Let ${\mathfrak A}_0$ be a $*$-subalgebra of
${\mathfrak A}$ and $\widetilde {\mathfrak A}_0$ be obtained by
adjoining a unit to ${\mathfrak A}_0$. Then ${\mathfrak A}_0$ is
said to be a {\em smooth subalgebra} of ${\mathfrak A}$ if the
following two conditions are satisfied:

\begin{enumerate}
\item ${\mathfrak A}_0$ is a dense $*$-subalgebra of ${\mathfrak A}$;

\item ${\mathfrak A}_0$ is stable under the holomorphic functional calculus,
that is, for any $a\in \widetilde {\mathfrak A}_0$ and for any
function $f$ that is holomorphic in a neighbourhood of the
spectrum of $a$ (thought of as an element in $\widetilde
{\mathfrak A}$) one has $f(a) \in \widetilde {\mathfrak A}_0$.
\end{enumerate}

Assume that ${\mathfrak A}_0$ is a dense $*$-subalgebra of
${\mathfrak A}$ such that ${\mathfrak A}_0$ is a Fr\'echet algebra
with a topology that is finer than that of ${\mathfrak A}$. A
necessary and sufficient condition for ${\mathfrak A}_0$ to be a
smooth subalgebra is given by the {\em spectral invariance}
condition cf. \cite[Lemma 1.2]{Schw}:

\begin{itemize}
\item $\widetilde {\mathfrak A}_0\cap GL(\widetilde {\mathfrak A}) =
GL(\widetilde {\mathfrak A}_0)$, where $GL(\widetilde{\mathfrak
A}_0)$ and $GL(\widetilde {\mathfrak A})$ denote the group of
invertibles in $\widetilde{\mathfrak A}_0$ and $\widetilde
{\mathfrak A}$ respectively.
\end{itemize}

\begin{rem}
This fact remains true in the case when ${\mathfrak A}$ is a
locally multiplicatively convex (i.e its topology is given by a
countable family of submultiplicative seminorms) Fr\'echet algebra
such that the group $GL(\widetilde {\mathfrak A})$ of invertibles
is open \cite[Lemma 1.2]{Schw}.
\end{rem}

One useful property of smooth subalgebras is the following. If
${\mathfrak A}_0$ is a smooth subalgebra of a $C^*$-algebra
${\mathfrak A}$, then the inclusion map ${\mathfrak A}_0\to
{\mathfrak A}$ induces an isomorphism in $K$-theory, \cite[Sect.
VI.3]{Co81}, \cite{Bost}. Another useful property of smooth
subalgebras is that sometimes there are interesting cyclic
cocycles on ${\mathfrak A}_0$ that do not extend to ${\mathfrak
A}$.

Let $\Gamma$ be a discrete group and $\sigma$ a multiplier on
$\Gamma$, $\ell$ denote the word length function on the group
$\Gamma$ with respect to a finite set of generators, i.e.
$\ell(\gamma) = d_\Gamma (\gamma, e)$ where $d_\Gamma$ denotes the
word metric. Let $\Delta$ denote the (unbounded) self-adjoint
operator on $\ell^2(\mathbb N)$ defined by $\Delta\delta_j = j
\delta_j$ for all $j\in \mathbb N$, and $D$ denote the (unbounded)
self-adjoint operator on $\ell^2(\Gamma)$ defined by
$D\delta_\gamma = \ell(\gamma) \delta_\gamma$ for all
$\gamma\in\Gamma $. Consider the unbounded derivation
$\widetilde\partial = {\rm ad}(D\otimes \id)$ of
$\B(\ell^2(\Gamma)\otimes \ell^2(\mathbb N))$. Recall that
$\widetilde\partial$ is a closed derivation with the domain ${\rm
Dom}\,(\widetilde\partial)$, which consists of all operators $T\in
\B(\ell^2(\Gamma)\otimes \ell^2(\mathbb N))$ such that $T$ maps
${\rm Dom}(D\otimes \id)$ into itself, and the operator
$\widetilde\partial(T)=(D\otimes \id)\circ T-T\circ (D\otimes
\id)$ defined initially on ${\rm Dom}(D\otimes \id)$ extends to a
bounded operator in $\ell^2(\Gamma)\otimes \ell^2(\mathbb N)$. Let
$$\B_\infty (\Gamma, \sigma) = \bigcap_{k\in \mathbb N} {\rm Dom}
(\widetilde\partial^k) \cap C^*_r(\Gamma, \bar\sigma)\otimes
{\mathcal K}$$ and $$ \B(\Gamma, \sigma) = \left\{T \in \B_\infty
(\Gamma, \sigma): \widetilde\partial^k(T)\circ (\id\otimes\Delta)
\;\;\;{\text {is bounded}}\ \forall k\in\N\right\}. $$ Then
$\B(\Gamma, \sigma) $ is a left ideal in $\B_\infty (\Gamma,
\sigma)$ - this follows from the observation that, since
$\widetilde\partial$ is a derivation,
\begin{equation}\label{e:ideal}
\widetilde\partial^k(T\circ S)\circ
(\id\otimes\Delta)=\sum_{m=0}^k\binom{k}{m}\widetilde
\partial^m(T)\circ
\widetilde\partial^{k-m}(S)\circ (\id\otimes\Delta).
\end{equation}

One has the following sufficient conditions for an operator $A\in
{\mathcal A}^L_{\ell^2(\N)}(\Gamma,\bar\sigma)$ to belong to the
algebra $\B(\Gamma,\sigma)$. For any $T\in \B(\ell^2(\N))$, let
$T_{ij}=(T(\delta_j),\delta_i), i,j\in \N,$ be the matrix elements
of $T$. Let $\mathcal R$ be the subalgebra in $\mathcal K$, which
consists of all compact operators in $\ell^2(\N)$, which are given
by rapidly decaying matrices, $$ \mathcal R = \Big\{T\in {\mathcal
K} : {\rm sup}\Big\{ i^k j^l |T_{ij}| : i,j\in \mathbb N\Big\} <
\infty, \;\;\;\forall \;k, l\in \mathbb N\Big\} . $$

\begin{lemma}\label{D2}
If $A\in {\mathcal A}^L_{\ell^2(\N)}(\Gamma,\bar\sigma)$,
$A=\sum_{\gamma\in\Gamma}T^L_\gamma\otimes A(\gamma)$ is such that
$A(\gamma) \in \mathcal R$ and also satisfies
\begin{equation}\label{D21}
\sum_{\gamma} \ell(\gamma)^k \| A(\gamma)\Delta\|<\infty,
\end{equation}
for all positive integers $k$, then $A \in  \B(\Gamma, \sigma)$.
\end{lemma}

\begin{proof} As in the proof of Lemma~\ref{D1}, let $K_1 \subset K_2 \subset
\cdots $ be a sequence of finite subsets of $\Gamma$ which is an
exhaustion of $\Gamma$, i.e. $\bigcup_{j\ge 1} K_j = \Gamma$. For
all $j \in \mathbb N$, define $A_j \in {\mathcal
A}^L_{\ell^2(\N)}(\Gamma,\bar\sigma)$ as
$A_j=\sum_{\gamma\in\Gamma}T^L_\gamma\otimes A_j(\gamma)$, where
$$ { {A_j}} (\gamma) = \left\{\begin{array}{l} {{A}} (\gamma) \;\;
{\rm if} \;\; \gamma\in K_j;\\[7pt]
         0 \;\;  {\rm otherwise.}
\end{array}\right.
$$ Then $A_j\in \C(\Gamma,\bar\sigma)\otimes {\mathcal R}$, and,
by Lemma~\ref{D1}, the sequence $A_j$ converges to $A$ in the norm
topology of ${\mathcal B}(\ell^2(\Gamma)\otimes \ell^2(\N))$.

Let $T\in \C(\Gamma,\bar\sigma)\otimes {\mathcal K}$. Since
$D\delta_e =0$, for any $v\in \ell^2(\N)$, one has
\begin{eqnarray*}
\widetilde\partial^k(T) (\delta_e\otimes v) &=& (D\otimes
\id)^k\circ T)(\delta_e\otimes v)\\ &=& (D\otimes \id)^k
\sum_{\gamma\in\Gamma}\delta_\gamma\otimes T(\gamma)v\\ &=&
\sum_{\gamma\in\Gamma}\delta_\gamma\otimes \ell(\gamma)^k
T(\gamma)v.
\end{eqnarray*}
Therefore, $T$ belongs to ${\rm Dom} (\widetilde\partial^k)$ for
any $k\in \mathbb N$, and
\begin{equation}\label{T}
[{\tilde\partial^k T}](\gamma) = \ell(\gamma)^k {T}(\gamma), \quad
\gamma\in \Gamma.
\end{equation}

As in the proof of Lemma~\ref{D1}, using Lemma~\ref{C},
(\ref{D21}) and (\ref{T}), one can establish that, for any $k\in
\N$, the sequences $\tilde\partial^k (A_j)$ and $\tilde\partial^k
(A_j) \circ (\id\otimes\Delta)$ converge in the norm topology of
${\mathcal B}(\ell^2(\Gamma)\otimes \ell^2(\N))$. This proves that
$A$ belongs to ${\rm Dom} (\widetilde\partial^k)$ for any $k\in
\mathbb N$, and the operator $\tilde\partial^k (A) \circ
(\id\otimes\Delta)$ is bounded in $\ell^2(\Gamma)\otimes
\ell^2(\N)$. Therefore $A\in \B(\Gamma, \sigma).$
\end{proof}

Following the arguments given in \cite{Co2}, III.5.$\gamma$, we
get

\begin{lemma}
The $*$-algebra $\B(\Gamma, \sigma)$ is a smooth subalgebra of
$C^*_r(\Gamma, \bar\sigma)\otimes {\mathcal K}$.
\end{lemma}
\begin{proof}
By Lemma~\ref{D2}, it follows that $\B(\Gamma, \sigma) $ contains
$\C(\Gamma,\bar\sigma)\otimes \mathcal R$. Therefore $\B(\Gamma,
\sigma) $ is a dense $*$-subalgebra of $C^*$-algebra
$C^*_r(\Gamma, \bar\sigma)\otimes {\mathcal K}$.

It is well-known (see, for instance, \cite[Theorem 1.2]{Ji92})
that $\B_\infty (\Gamma, \sigma)$ is a smooth subalgebra of
$C^*_r(\Gamma, \bar\sigma)\otimes {\mathcal K}$. Since $\B(\Gamma,
\sigma)$ is a left ideal in $\B_\infty (\Gamma, \sigma)$, the
proof is completed by the following simple algebraic fact: if $A$
is a spectral invariant subalgebra of an algebra $B$ and $I$ is a
left ideal in $A$, then $I$ is spectral invariant in $B$.
\end{proof}

For any $k\in \N$ and $f\in C^*_r(\Gamma,\bar\sigma)$, put
$$\nu_k(f) = (\sum_{\gamma\in \Gamma} (1+\ell(\gamma))^{2k}
|f(\gamma)|^2)^{\frac{1}{2}}, \quad k\in\N.$$ We clearly have that
$\nu_k(f)<\infty$ for any $k\in \N$ and
$f\in\C(\Gamma,\bar\sigma)$.

Consider any element $A$ in $C^*_r(\Gamma,\bar\sigma)
\otimes\mathcal K$, $A=\sum_{\gamma\in\Gamma}T^L_\gamma\otimes
A(\gamma)$. Let $(A_{ij}(\gamma))$ be the matrix, corresponding to
$A(\ga)$. Then $A_{ij} \in C^*_r(\Gamma, \bar\sigma)$ for any
$i\in\N$ and $j\in\N$. Put $$N_k(A) = ( \sum_{i,j}
\nu_k(A_{ij})^2)^{\frac{1}{2}}, \quad k\in \N.$$

\begin{lemma}
For all $A\in  \B(\Gamma, \sigma)$ and $k\in \mathbb N$, we have
$N_k(A)<\infty$ .
\end{lemma}
\begin{proof}
One can be easily seen that, for any $A\in
C^*_r(\Gamma,\bar\sigma) \otimes\mathcal K$,
\[
A(\delta_e\otimes\delta_j)= \sum_{i\in \mathbb N, \ga\in \Gamma}
A_{ij}(\ga)\delta_\ga\otimes \delta_i.
\]
Using (\ref{T}), one has
\begin{eqnarray*}
\widetilde\partial^k(A)\circ (\id\otimes\Delta)
(\delta_e\otimes\delta_j) &=& j\, \widetilde\partial^k(A)
(\delta_e\otimes \delta_j)  \\ &=&  j\, \sum_{i\in \mathbb N,
\ga\in \Gamma} \ell(\ga)^k A_{ij}(\ga)\delta_\ga\otimes \delta_i.
\end{eqnarray*}
Therefore for $A\in  \B(\Gamma, \sigma)$ and $k\in \mathbb N$,
there is a positive constant $C$ such that $$ \sum_{i\in \mathbb
N, \ga\in \Gamma} \ell(\ga)^{2k}  |A_{ij}(\ga)|^2 <Cj^{-2}. $$
Hence, $$ \sum_{i, j\in \mathbb N, \ga\in \Gamma} \ell(\ga)^{2k}
|A_{ij}(\ga)|^2 <C, $$ that completes the proof.
\end{proof}

\subsection{Group cocycles and cyclic cocycles}
The cyclic cocycles that we consider arise from normalised group
cocycles on $\Gamma$. Recall (see, for instance, \cite{Gui}) that
a (homogeneous) group $k$-cocycle is a map $h:\Gamma^{k+1} \to \C$
satisfying the identities
\begin{gather*}
h(\ga \ga_0, \ldots ,\ga \ga_k) = h(\ga_0, \ldots ,\ga_k);\\
\sum_{i=0}^{k+1} (-1)^i h(\ga_0, \ldots,
\ga_{i-1},\ga_{i+1},\ldots, \ga_{k+1}) = 0.
\end{gather*}
Then an (inhomogeneous) group $k$-cocycle $c\in Z^k(\Gamma,
\mathbb C)$ that is associated to such an $h$ is given by $$
c(\ga_1, \ldots, \ga_k) = h (e, \ga_1, \ga_1 \ga_2, \ldots,
\ga_1\ldots \ga_k). $$ It can be easily checked that $c$ satisfies
the following identity
\begin{equation}\label{e:c}
c(\ga_1, \ga_2,\ldots, \ga_k)
+\sum_{i=0}^{k-1}(-1)^{i+1}c(\ga_0,\ldots, \ga_{i-1},
\ga_i\ga_{i+1}, \ga_{i+2},\ldots,\ga_k) + (-1)^{k+1}
c(\ga_0,\ga_1,\ldots,\ga_{k-1})=0.
\end{equation}

A group $k$-cocycle is said to be normalised (in the sense of
Connes), if $c(\ga_1,\ga_2,\ldots,\ga_k)$ is zero if either $\ga_i
= e$ for some $i$ or if $\ga_1\ldots \ga_k = e$.

Recall that a cyclic $k$-cocycle on an algebra $A$ is
a $k+1$-linear functional $\phi$ on $A$, satisfying the following
identities
\begin{gather}
\phi(f_k,f_0,\ldots,f_{k-1})=(-1)^k\phi(f_0,f_1,\ldots,f_k),
\label{e;cyclic1}\\
\begin{split}
b\phi(f_0,f_1,\ldots,f_{k+1}) &\equiv\sum_{0\leq j\leq k}(-1)^j
\phi(f_0,\ldots,f_jf_{j+1},\ldots,f_{k+1})+(-1)^{k+1}
\phi(f_{k+1}f_0,f_1,\ldots,f_k)\\ &=0 \label{e:cyclic2},
\end{split}
\end{gather}
where $f_0, f_1,\ldots, f_{k+1}\in A$.

Given a normalised group cocycle $c\in Z^k(\Gamma, \mathbb C)$,
$k=0,\ldots, \dim M$, we define a cyclic $k$-cocycle $\tau_c$ on
the twisted group ring $\mathbb C(\Gamma, \bar\sigma)$, which is
given by
\begin{equation}
\label{e:tauc} \tau_c(a_0\delta_{\ga_0}, \ldots,a_k\delta_{\ga_k})
= \left\{\begin{array}{l} a_0\ldots a_k c(\ga_1,\ldots, \ga_k)
\tr_\Gamma( \delta_{\ga_0}\ast\delta_{\ga_1}\ast
\ldots\ast\delta_{\ga_k}) \quad\text{if} \,\, \ga_0 \ldots \ga_k
=e;\\
\\
0 \qquad \text{otherwise,}
\end{array}\right.
\end{equation}
where $a_j \in \C$ for $j=0,1,\ldots k$.

\begin{lemma}
For any normalised group cocycle $c\in Z^k(\Gamma, \mathbb C)$,
$k=0,\ldots, \dim M$, the functional $\tau_c$ defined by the
formula (\ref{e:tauc}) is a cyclic cocycle on
$\C(\Gamma,\bar\sigma)$.
\end{lemma}

\begin{proof}
It is clearly sufficient to check the identity (\ref{e;cyclic1})
in the case when $f_j=a_j\delta_{\ga_j}, j=0,1,\ldots,k,$ with
$a_j\in \C$ and $\ga_0\ga_1\ldots \ga_k=e$. Then, since $c$ is
normalised, (\ref{e:c}) implies that $$ c(\ga_1, \ga_2,\ldots,
\ga_k)+(-1)^{k+1}c(\ga_0,\ga_1,\ldots,\ga_{k-1})=0, $$ from where
(\ref{e;cyclic1}) follows immediately.

To prove the identity (\ref{e:cyclic2}), we again assume that
$f_j=a_j\delta_{\ga_j}, j=0,1,\ldots,k+1,$ with $a_j\in \C$ and
$\ga_0\ga_1\ldots \ga_{k+1}=e$. Then we have
\begin{align*}
\tau_c(f_0f_1,f_2,\ldots,f_{k+1})& =\tau_c(a_0\delta_{\ga_0}\ast
a_1\delta_{\ga_1},a_2\delta_{\ga_2},\ldots,
a_{k+1}\delta_{\ga_{k+1}}) \\ &=\sigma (\ga_0,\ga_1)
\tau_c(a_0a_1\delta_{\ga_0\ga_1},a_2\delta_{\ga_2},\ldots,
a_{k+1}\delta_{\ga_{k+1}})\\ &=\sigma (\ga_0,\ga_1) a_0\ldots
a_{k+1} c(\ga_2,\ldots,\ga_{k+1}) \tr_\Gamma(
\delta_{\ga_0\ga_1}\ast\delta_{\ga_2}\ast
\ldots\ast\delta_{\ga_{k+1}})\\ &=a_0\ldots a_{k+1}
c(\ga_2,\ldots,\ga_{k+1}) \tr_\Gamma(
\delta_{\ga_0}\ast\delta_{\ga_1}\ast\delta_{\ga_2}\ast
\ldots\ast\delta_{\ga_{k+1}}).
\end{align*}
Similarly, we get for $j=1,2,\ldots,k$ $$
\tau_c(f_0,\ldots,f_jf_{j+1},\ldots,f_{k+1}) =a_0\ldots a_{k+1}
c(\ga_1,\ldots,\ga_j\ga_{j+1},\ldots,\ga_{k+1}) \tr_\Gamma(
\delta_{\ga_0}\ast\delta_{\ga_1}\ast\delta_{\ga_2}\ast
\ldots\ast\delta_{\ga_{k+1}}), $$ and $$
\tau_c(f_{k+1}f_0,f_1,\ldots,f_k)= a_0\ldots a_{k+1}
c(\ga_1,\ldots,\ga_{k})
\tr_\Gamma(\delta_{\ga_0}\ast\delta_{\ga_1}\ast\delta_{\ga_2}\ast
\ldots\ast\delta_{\ga_{k+1}}). $$ Using these identities and
(\ref{e:c}), we obtain
\begin{eqnarray*}
b\tau_c(f_0,f_1,\ldots,f_{k+1}) &=&a_0\ldots a_{k+1}\tr_\Gamma(
\delta_{\ga_0}\ast\delta_{\ga_1}\ast\delta_{\ga_2}\ast
\ldots\ast\delta_{\ga_{k+1}})\\ & &(c(\ga_2,\ldots,\ga_{k+1})+
\sum_{j=1}^k
(-1)^jc(\ga_1,\ldots,\ga_j\ga_{j+1},\ldots,\ga_{k+1})\\ &
&\mbox{}+(-1)^{k+1} c(\ga_1,\ldots,\ga_{k}))=0,
\end{eqnarray*}
that completes the proof.
\end{proof}

Of particular interest is the case $k=2$, when the formula
(\ref{e:tauc}) reduces to $$ \tau_c(a_0\delta_{\ga_0},
a_1\delta_{\ga_1}, a_2\delta_{\ga_2})  = \left\{\begin{array}{l}
        a_0 a_1 a_2 c(\ga_1, \ga_2)
\sigma(\ga_1, \ga_2) \sigma(\ga_2^{-1}, \ga_2)\quad \text{if} \,\,
\ga_0 \ga_1 \ga_2 =e ;\\
\\
0 \qquad \text{otherwise.}
\end{array}\right.
$$

For any normalised group cocycle $c\in Z^k(\Gamma, \mathbb C)$,
$k=0,\ldots, \dim M$, one can define a cyclic $k$-cocycle
$\tau_c\#\Tr$ on $\C(\Gamma, \bar\sigma)\otimes {\mathcal R}$ by
the formula
\begin{equation}\label{e:cup}
\tau_c\#\Tr(f_0,f_1,\ldots,f_k)=\sum_{\ga_0\ga_1\ldots \ga_k=e}
\Tr(f_0(\ga_0)f_1(\ga_1)\ldots f_k(\ga_k))c(\ga_1,\ldots,\ga_k)
\tr_\Gamma(\delta_{\ga_0}\ast\delta_{\ga_1}\ast\ldots\ast\delta_{\ga_k}),
\end{equation}
where $f_j=\sum_{\ga\in\Gamma}\delta_\ga\otimes
f_j(\ga)\in\C(\Gamma,\bar\sigma)\otimes {\mathcal R},
j=0,1,\ldots,k.$

Recall that a normalised $k$-cocycle $c$ is said to be
polynomially bounded if there are a positive constant $C$ and
$a_i\in \mathbb N$ for all $i=1, \ldots, k$, such that
\begin{equation}\label{polynomial}
|c(\ga_1, \ga_2, \ldots, \ga_k)|\le C (1+\ell(\ga_1))^{a_1}
(1+\ell(\ga_2))^{a_2} \ldots (1+\ell(\ga_k))^{a_k}.
\end{equation}
For example, it is easy to see that any normalised group
$1$-cocycle is polynomially bounded. Groups $\Gamma$ that are
virtually nilpotent or that are word hyperbolic have the property
that every group cohomology class has  a representative cocycle
that is polynomially bounded, cf. \cite{Gr, Gr2}.

Recall that a discrete group $\Gamma$ has property (RD), if the
Haagerup inequality holds for $\Gamma$, that is, there exist $N\in
\N$ and $C'>0$ such that
\begin{equation}\label{Haagerup}
\|f\|_{C^*_r(\Gamma)}\leq C'\nu_N(f), \quad f\in \C\Gamma,
\end{equation}
or, equivalently,
\[
\|f*_0u\|\leq C' \nu_N(f)\|u\|, \quad f\in \C\Gamma, \quad u\in
\ell^2(\Gamma),
\]
where $\|\cdot\|$ denotes the norm in $\ell^2(\Gamma)$:
$\|u\|=\nu_0(u)$, and $*_0$ denotes the usual (untwisted)
convolution in $\C\Gamma$:
\[
f*_0u(\ga)=\sum_{\ga_1\ga_2=\ga}f(\ga_1)u(\ga_2), \quad \gamma\in
\Gamma.
\]
Examples of groups, having property (RD), are virtually
nilpotent groups, hyperbolic groups and products of hyperbolic
groups, cocompact lattices in $\rm{SL}_3(\R), \rm{SL}_3(\C),
\rm{SL}_3(\mathbb H)$.

One shows exactly as in \cite{Co2} section III.5.$\gamma$ using
the Haagerup inequality (\ref{Haagerup}), that if $\Gamma$ has
property (RD), then for each polynomially bounded normalised group
cocycle $c$, the associated cyclic cocycle $\tau_c$ on ${\mathbb
C}(\Gamma, \bar\sigma)$, is continuous for the norm $\nu_K$, for
$K$ sufficiently large, and the tensor product cocycle
$\tau_c\#\Tr$ defined on $\C(\Gamma,\bar\sigma)\otimes {\mathcal
R}$ by the formula (\ref{e:cup}) extends by continuity to
$\B(\Gamma, \sigma)$.

\begin{lemma}\label{taucforrd}
Let a group $\Gamma$ have property (RD) and a normalised group
$k$-cocycle $c$ be polynomially bounded. Then the associated
cyclic cocycle $\tau_c$ on ${\mathbb C}(\Gamma, \bar\sigma)$, is
continuous for the norm $\nu_K$, for $K$ sufficiently large.
\end{lemma}

\begin{proof}
For all $f_0, f_1, \ldots f_k \in \C(\Gamma, \bar\sigma)$, one
has, $$ \tau_c(f_0, f_1, \ldots f_k) = \sum_{\ga_0\ga_1\ldots
\ga_k=e} f_0(\ga_0)\ldots f_k(\ga_k) c(\ga_1, \ldots \ga_k)
\tr_\Gamma(\delta_{\ga_0} *\delta_{\ga_1}\*\ldots \delta_{\ga_k}).
$$

Using the Haagerup inequality (\ref{Haagerup}), we get
\begin{align*}
\left|  \tau_c(f_0, f_1, \ldots f_k) \right| & =  \left|
\tau_c(f_1, f_2, \ldots f_k, f_0) \right|\\ &\le  C
\sum_{\ga_0\ga_1\ldots \ga_k=e} |f_1(\ga_1)|\ldots |f_k(\ga_k)|
|f_0(\ga_0)| (1+\ell(\ga_1))^{a_1} (1+\ell(\ga_2))^{a_2} \ldots
(1+\ell(\ga_k))^{a_k}
     \\
&=C  |(1+\ell)^{a_1} f_1|*_0\cdots *_0|(1+\ell)^{a_k}
f_k|*_0|f_0|(e)
\\ &\le C (C')^n \|f_0\| \prod_{j=1}^k \nu_N(|(1+\ell)^{a_j}
f_j|)\\ &= C (C')^n \nu_0(f_0) \prod_{j=1}^k \nu_{N+a_j}(|f_j|)\\
&= C (C')^n \nu_0(f_0) \prod_{j=1}^k \nu_{N+a_j}(f_j),
\end{align*}
that proves the desired continuity property.
\end{proof}

\begin{lemma}\label{RD}
If $\Gamma$ has property (RD), and given a polynomially bounded
normalised group $k$-cocycle $c$, a cyclic $k$-cocycle
$\tau_c\#\Tr$ on $\C(\Gamma, \bar\sigma)\otimes {\mathcal R}$
given by the formula (\ref{e:cup}) extends by continuity to
$\B(\Gamma, \sigma)$.
\end{lemma}

\begin{proof}
The proof is a word-by-word repetition of the proof of Lemma~6.4
in \cite{CM90}. Take any $f_l=\sum_{\ga\in\Gamma}\delta_\ga\otimes
f_l(\ga)\in\C(\Gamma,\bar\sigma)\otimes {\mathcal R},
l=0,1,\ldots,k.$ Represent any $f_l(\ga)\in \mathcal R$ by a
matrix $(f^l_{ij}(\ga))$. Then, for any $l=0,1,\ldots,k,$ $i\in\N$
and $j\in\N$, we have $f^l_{ij}\in \C(\Gamma,\bar\sigma)$, and one
can see that $$
\tau_c\#\Tr(f_0,f_1,\ldots,f_k)=\sum_{i_0,i_1,\ldots,i_k}
\tau_c(f^0_{i_0i_1},f^1_{i_1i_2},\ldots ,f^k_{i_ki_0}). $$ So,
using Lemma~\ref{taucforrd}, we have
\begin{eqnarray*}
|\tau_c\#\Tr(f_0,f_1,\ldots,f_k)|&\leq &\sum_{i_0,i_1,\ldots,i_k}
|\tau_c(f^0_{i_0i_1},f^1_{i_1i_2},\ldots ,f^k_{i_ki_0})|\\ &\leq &
C \sum_{i_0,i_1,\ldots,i_k}
\nu_K(f^0_{i_0i_1})\nu_K(f^1_{i_1i_2})\ldots \nu_K(f^k_{i_ki_0})
\end{eqnarray*}
with some natural $K$. Then we use the following inequality $$
\sum_{i_0,i_1,\ldots,i_k} \alpha^0_{i_0i_1}\alpha^1_{i_1i_2}\ldots
\alpha^k_{i_ki_0} \leq
\prod_{l=0}^k\left(\sum_{i,j}(\alpha_{ij}^l)^2\right)^{1/2}, $$
which holds for any $k\geq 1$, that gives us the estimate $$
|\tau_c\#\Tr(f_0,f_1,\ldots,f_k)|\leq CN_K(f_0)N_K(f_1)\ldots
N_K(f_k), $$ and concludes the proof.
\end{proof}

\subsection{A vanishing lemma for pairing with cyclic cocycles}
Let a group $\Gamma$ be a discrete group, $c$ be a normalised
group $k$-cocycle on $\Gamma$ ($k$ even), and $\tau_c$ be the
associated cyclic cocycle on $\C(\Gamma, \bar\sigma)$. By the
pairing theory of \cite{Co} we get an additive map
\[
          [\tau_c\#\Tr] : K_0(\C(\Gamma, \bar\sigma)\otimes{\mathcal R})
          \to\mathbb{R}.
\]
Explicitly, $[\tau_c\#\Tr]([e]-[f]) = \widetilde\tau_c(e,\cdots,e)
- \widetilde\tau_c(f,\cdots,f)$, where $e,f$ are idempotent
matrices with entries in $\C(\Gamma,
\bar\sigma)\otimes\widetilde{\mathcal R}$, which is the unital
algebra obtained by adding the identity to $\C(\Gamma,
\bar\sigma)\otimes{\mathcal R}$, and $\widetilde\tau_c$ denotes
the canonical extension of $\tau_c\#\Tr$ to $\C(\Gamma,
\bar\sigma)\otimes\widetilde{\mathcal R}\otimes M_N(\C)$ defined
as follows; $$
        \widetilde\tau_c(f_0\otimes R_0, \ldots, f_k\otimes R_k)
= \tr(R_0\ldots R_k) \; \tau_c\#\Tr(f_0, \ldots, f_k), $$ where
$f_j\in \C(\Gamma, \bar\sigma)\otimes\widetilde{\mathcal R},
R_j\in M_N(\C), j=0,1,\ldots,k.$

Recall that $\Tr_{\Gamma}$ denotes the tensor product of the
canonical finite trace $\tr_\Gamma$ on $\C(\Gamma, \bar\sigma)$
and the standard trace $\Tr$ on $\mathcal R$.

\begin{lemma}\label{l:vanishing}
Let $c$ be a normalised $k$-cocycle on a discrete group $\Gamma$
($k$ even) and $\tau_c\#\Tr$ be the associated cyclic $k$-cocycle
on the group algebra $\C(\Gamma, \bar\sigma)\otimes{\mathcal R}$.
Let $I\otimes P$ be the tensor product of the identity in
$\C(\Gamma, \bar\sigma)$ and a projection $P$ in $\mathcal R$.
Then we have, $$
\begin{array}{rcl}
\Tr_{\Gamma}(I\otimes P) &=&  {\rm rank}(P);\\[+7pt]
\tau_c\#\Tr(I\otimes P,\ldots, I\otimes P) &=& 0\qquad {\text{for
$k>0$}}.
\end{array}
$$
\end{lemma}

\begin{proof} Observe that $I\otimes P = \delta_e\otimes P$.
The statement is trivially true when $k=0$. For $k>0$, we have $$
\begin{array}{rcl}
\tau_c\#\Tr(I\otimes P,\ldots, I\otimes P) &=&
\tau_c\#\Tr(\delta_e\otimes P,\ldots, \delta_e\otimes P)
\\[+7pt]
&=& {\rm rank}(P) \tau_c(\delta_e,\ldots,\delta_e)
\\[+7pt]
& = &  {\rm rank}(P)  c(e,\ldots,e) \tr_\Gamma( \delta_{e}\ast\delta_{e}\ast
\ldots\ast\delta_{e}).
\end{array}
$$
Since $c$ is a normalised group cocycle, $c(e,\ldots,e)=0$, and the result
follows.
\end{proof}

\section{Semiclassical vanishing theorems for spectral projections}
\label{vanishing}
\subsection{General results on equivalence of projections in smooth
subalgebras} In the setting of Section~\ref{s:abstract-equiv},
suppose, in addition, that there is given a smooth subalgebra
${\mathfrak A}_0$ in the $C^*$-algebra ${\mathfrak A}$ such that:
\begin{itemize}
\item for any $t>0$, the operators $e^{-tA_l}$ belong to
$\pi_l({\mathfrak A}_0)$, $l=1,2$.
\end{itemize}

Consider an interval $(a_1,b_1)$, which does not intersect with
the spectrum of $A_1$. Let $a_2, b_2$ be given by the formulas
\eqref{e:a2} and \eqref{e:b2}. Suppose that
$\alpha_1>a_1+\gamma_1$, $\alpha_2>b_2+\gamma_2$ and $b_2>a_2$. By
Theorem~\ref{t:equivalence}, the interval $(a_2,b_2)$ does not
intersect with the spectrum of $A_2$.

Note that $E_1(\lambda_1) = \chi_{[e^{-t\lambda_1},
\infty)}\left(e^{-tA_1}\right)$. Since  $\lambda_1\not\in {\rm
spec}(A_1)$, by the Riesz formula one has, $$ E_1(\lambda_1) =
\frac{1}{2\pi i} \oint_C (\lambda - e^{-t A_1})^{-1} d\lambda, $$
where $C$ is a contour intersecting the real axis at
$e^{-t\lambda_1}$ and at some large positive number not in the
spectrum of $e^{-t A_1}$. It follows that $E_1(\lambda_1)$ is a
holomorphic function of $e^{-t A_1}$, and therefore one has
$E_1(\lambda_1)\in\pi_1({\mathfrak A}_0)$. Similarly, for any
$\lambda_2\in (a_2,b_2)$, the spectral projection $E_2(\lambda_2)$
belongs to $\pi_2({\mathfrak A}_0)$.

\begin{thm}\label{t:equivalence1}
The projections ${\mathcal V}_1E_1(\lambda_1){\mathcal V}_1^{-1}$
and ${\mathcal V}_2E_2(\lambda_2){\mathcal V}_2^{-1}$ are
Murray-von Neumann equivalent in $\pi({\mathfrak A}_0)$ for any
$\lambda_1\in (a_1,b_1)$ and $\lambda_2\in (a_2,b_2)$.
\end{thm}

\begin{proof}
Consider a bounded operator $T={\mathcal V}_2E_2(\lambda_2)i_2 J
p_1 E_1(\lambda_1) {\mathcal V}_1^{-1}$ in $\chh$. As shown in the
proof of Theorem~\ref{t:equivalence}, the operator $T$ belongs to
$\pi({\mathfrak A})$ and is invertible as an operator from
${\mathcal V}_1 (\Im E_1(\lambda_1))$ to ${\mathcal V}_2 (\Im
E_2(\lambda_2))$. Since $\pi({\mathfrak A}_0)$ is dense in
$\pi({\mathfrak A})$, there exists an operator $T_1$ in
$\pi({\mathfrak A}_0)$ such that the operator $P={\mathcal V}_2
E_2(\lambda_2){\mathcal V}^{-1}_2 T_1{\mathcal V}_1
E_1(\lambda_1){\mathcal V}^{-1}_1$ is invertible as an operator
from ${\mathcal V}_1 (\Im E_1(\lambda_1))$ to ${\mathcal V}_2 (\Im
E_2(\lambda_2))$. Then the operator $P$ as an operator in $\chh$
belongs to $\pi({\mathfrak A}_0)$, and its image $\Im P={\mathcal
V}_2 (\Im E_2(\lambda_2))$ is closed. The desired statement
follows from the following lemma.

\begin{lemma}\label{l:polar}
Let ${\mathfrak A}$ be a $C^*$-algebra, $\chh$ a Hilbert space
equipped with a faithful $\ast$-representation of ${\mathfrak A}$,
$\pi: {\mathfrak A}\to \ck(\chh)$, ${\mathfrak A}_0$ a smooth
subalgebra in ${\mathfrak A}$. If a bounded operator $P$ in $\chh$
belongs to $\pi({\mathfrak A}_0)$ and has closed image, and $P=US$
is its polar decomposition, then $U, S\in \pi({\mathfrak A}_0)$.
\end{lemma}

\begin{proof}
As shown in the proof of Lemma~\ref{l:polar1}, zero is an isolated
point in the spectrum of $P^*P$. Since ${\mathfrak A}_0$ is stable
under holomorphic functional calculus, the operator
$S=\sqrt{P^*P}$ is in ${\mathfrak A}_0$ and has zero as an
isolated point in its spectrum. Furthermore, the function $f$
introduced in the proof of Lemma~\ref{l:polar1} extends to a
holomorphic function in a neighborhood of the spectrum of $S$,
that implies that $S^{(-1)}=f(S)$ and $U=PS^{(-1)}$ are also in
$\pi({\mathfrak A}_0)$.
\end{proof}

Applying Lemma~\ref{l:polar} to the operator $P$ as above, we
obtain an isometry $U\in\pi({\mathfrak A}_0)$, which gives the
desired equivalence of the projections ${\mathcal
V}_1E_1(\lambda_1){\mathcal V}_1^{-1}$ and ${\mathcal
V}_2E_2(\lambda_2){\mathcal V}_2^{-1}$.
\end{proof}

\subsection{Proof of the semiclassical vanishing theorem of the higher
traces of spectral projections} In this Section, we prove the
second part of Theorem~\ref{main1} and Corollary~\ref{main2}. For
this, we apply Theorem~\ref{t:equivalence1} in the setting of
Section~\ref{s:model1}. We will use the notation of this section.
Here we make a particular choice of the unitary isomorphisms
$V_1:L^2({\mathbb R}^n,\C^k)^N\to \ell^2(\N)$ and
$V_2:L^2(\F,\tE|_\F)\to \ell^2(\N)$. Namely, we define $V_1$,
using the spectral decomposition for the model operator $K=K(1)$
in $L^2({\mathbb R}^n,\C^k)^N$. More precisely, suppose that
$\lambda_1\leq \lambda_2\leq \lambda_3\leq \ldots $ are the
eigenvalues of the operator $K$ (counting with multiplicities),
and $\phi_j$ the corresponding eigenfunctions, which form a
complete orthonormal system in $L^2({\mathbb R}^n,\C^k)^N$. Then
$V_1$ is defined as $V_1\phi_j=\delta_j, j\in \N$. Similarly, take
any second order self-adjoint $(\Gamma,\sigma)$-invariant elliptic
differential operator $P$ in $L^2(\M,\tE)$. Define the unitary
isomorphism $V_2:L^2(\F,\tE|_\F)\to \ell^2(\N)$, using the
orthonormal basis of eigenvectors for the operator $P$ in
$L^2(\F,\tE|_\F)\cong L^2(M,{\mathcal E})$ with the
$(\Gamma,\sigma)$-periodic boundary conditions as above.

Take a smooth subalgebra ${\mathfrak A}_0$ of the $C^*$-algebra
${\mathfrak A}=C^*_r(\Gamma, \bar\sigma)\otimes {\mathcal K}$ to
be $\B(\Gamma, \sigma)$.

\begin{lemma}\label{F1}
For any $t>0$, the operator ${\mathcal V}_1e^{-tA_1}{\mathcal
V}^*_1$ belongs to $\B(\Gamma, \sigma) \subset C^*_r(\Gamma,
\bar\sigma) \otimes\K$.
\end{lemma}

\begin{proof}
Since ${\mathcal V}_1e^{-tA_1}{\mathcal V}^*_1=\id\otimes
V_1e^{-tK(\mu)}V^*_1$, it suffices to prove that the operator
$V_1e^{-tK(\mu)}V^*_1\Delta$ is bounded in $\ell^2(\N)$. By
definition of $V_1$, the operator $V_1 K V^*_1$ in $\ell^2(\N)$ is
given by
\[
V_1 K V^*_1\delta_j=\lambda_j\delta_j, \quad j\in \N.
\]
It is well-known that $\lambda_j \sim Cj^{1/n}, j\to \infty $.
Therefore, the operator $K^{-n} V^*_1\Delta V_1$ is bounded in
$L^2({\mathbb R}^n,\C^k)^N$. Since $V_1e^{-tK(\mu)}V^*_1\Delta =
V_1e^{-tK(\mu)}K^n K^{-n}V^*_1\Delta V_1 V^*_1$, and the operator
$e^{-tK(\mu)}K^n$ is bounded in $L^2({\mathbb R}^n,\C^k)^N$, this
immediately completes the proof.
\end{proof}

\begin{lemma}\label{F}
For all $t>0$, ${\mathcal V}_2e^{-tH_{}(\mu)}{\mathcal V}^*_2 \in
\B(\Gamma, \sigma) \subset C^*_r(\Gamma, \bar\sigma) \otimes\K$.
\end{lemma}

\begin{proof} First, recall the following well-known properties of the heat
operator $e^{-tH_{}(\mu)}$, cf. \cite{Greiner,Ko}. Let $d$ denote
the Riemannian distance function on $\widetilde M$.

\begin{lemma}{\label E}
The Schwartz kernel $k(t,x,y)$ of the heat operator
$e^{-tH_{}(\mu)}$ is smooth for all $t>0$. Moreover, for any $t>0$
and for any $(\Gamma,\sigma)$-invariant differential operators
$A=a(x,D_x)$ and $B=b(x,D_x)$ in $C^\infty(\M, \tE)$ there are
positive constants $C_1, C_2$ depending on $\mu$ such that the
following off-diagonal estimate holds $$
|a(x,D_x)b(y,D_y)k(t,x,y)| \le C_1 e^{-C_2 d(x,y)^2},\quad x\in
\M,\quad y\in \M. $$
\end{lemma}

We have ${\mathcal V}_2e^{-tH_{}(\mu)}{\mathcal V}^*_2 \in
{\mathcal A}^L_{\ell^2(\N)}(\Gamma,\bar\sigma)$, so we can write
the operator ${\mathcal V}_2e^{-tH_{}(\mu)}{\mathcal V}^*_2$ as
\[
{\mathcal V}_2e^{-tH_{}(\mu)}{\mathcal V}^*_2
=\sum_{\gamma\in\Gamma}T^L_\gamma\otimes h_{t,\mu}(\gamma)
\]
with some $h_{t,\mu}(\gamma)\in\B(\ell^2(\N))$. We will identify
the space $L^2(\F,\tE|_\F)$ with the subspace in $L^2(\M, \tE)$,
which consists of sections from $L^2(\M, \tE)$, vanishing outside
of $\F$. As above, $i: \mathcal{F} \to \M$ denote the inclusion
map.

\begin{lemma}\label{F2}
The operator $V_2^*h_{t,\mu}(\gamma)V_2$ is given by the
restriction of the operator $i^*T_\gamma e^{-tH_{}(\mu)}$ to
$L^2(\F,\tE|_\F)$.
\end{lemma}

\begin{proof}
Recall that a unitary operator ${\mathcal V}_2 : \mathfrak{H} \to
{\mathcal H}$ is defined as ${\mathcal V}_2=(\id\otimes V_2)\circ
{\bf U}$, where ${\bf U}$ is the $(\Gamma, \sigma)$-equivariant
isometry \eqref{isom}. Therefore, we have
\[
{\bf U}e^{-tH_{}(\mu)}{\bf U}^*
=\sum_{\gamma\in\Gamma}T^L_\gamma\otimes
V_2^*h_{t,\mu}(\gamma)V_2.
\]
Let $\phi\in L^2(\F,\tE|_\F)$. By definition of $\bf U$, it easily
follows that ${\bf U}^*(\delta_e\otimes \phi)\in L^2(\M,\tE)$
coincides with $\phi$. Therefore, by \eqref{isom}, we get
\[
{\bf U}e^{-tH_{}(\mu)}{\bf U}^*(\delta_e\otimes
\phi)=\sum_{\gamma\in\Gamma}\delta_\gamma\otimes i^*T_\gamma
e^{-tH_{}(\mu)}\phi,
\]
that immediately completes the proof.
\end{proof}

As in the proof of Lemma~\ref{F1}, using the Weyl asymptotic
formula $\alpha_j\sim Cj^{2/n}, j\to \infty,$ for the eigenvalues
$\alpha_j$ of the operator $P$, one can show that the operator
$P^{-n/2}V_2^*\Delta V_2$ is bounded in $L^2(\F,\tE|_\F)$. By this
fact and Lemma~\ref{F2}, it follows that
\begin{equation}\label{e:3}
\|h_{t,\mu}(\gamma)\Delta\| =
\|V_2^*h_{t,\mu}(\gamma)V_2V_2^*\Delta V_2\| =\|i^*T_\gamma
e^{-tH_{}(\mu)}P^{m}P^{-m}V_2^*\Delta V_2\| \leq C_3\|i^*T_\gamma
e^{-tH_{}(\mu)}P^{m}\|
\end{equation}
for any natural $m>n/2$ with some positive constant $C_3$.

It is well known that
\begin{equation}\label{e:d}
\ell(\gamma) \le C_4 (\inf_{x,y\in\F} d (\gamma x, y) +1)
\end{equation}
for some positive constant $C_4$. From (\ref{e:d}) and
Lemma~\ref{E}, we get
\begin{equation}\label{e:2}
\|i^*T_\gamma e^{-tH_{}(\mu)}P^{m}\|\le C_5 e^{-C_6
\ell(\gamma)^2}.
\end{equation}
Observe that one has the estimate
\begin{equation}\label{e:d1}
  \# \left\{\gamma\in\Gamma\; |\;
\ell(\gamma) \le R\right\} \le C_7 e^{C_8 R}
\end{equation}
for some positive constants $C_7, C_8$, since the growth rate of
the volume of balls in $\Gamma$ is at most exponential. By
(\ref{e:3}), (\ref{e:2}) and (\ref{e:d1}), it follows that $$
\sum_\gamma \ell(\gamma)^k \|h_{t,\mu}(\gamma)\Delta\| <\infty $$
for all positive integers $k$. By Lemma \ref{D2}, this implies
that $\;{\mathcal V}_2e^{-tH_{}(\mu)}{\mathcal V}^*_2 \in
\B(\Gamma, \sigma) \subset C^*_r(\Gamma, \bar\sigma) \otimes\K$
for all $t>0$.
\end{proof}

\begin{rem}
Since the Schwartz kernel of $e^{-tH_{}(\mu)}$ is smooth $\forall
\gamma\in\Gamma$ by Lemma~\ref{E}, it follows from the proof of
Lemma~\ref{F} that $h_{t,\mu}(\gamma) \in \mathcal R$ $\forall
\gamma\in\Gamma$ (cf. also Lemma 5 in III.4.$\beta$ in
\cite{Co2}).
\end{rem}

By Lemmas~\ref{F1} and~\ref{F}, it follows that, for any $t>0$,
the operators $e^{-tA_l}$ belong to $\pi_l({\mathfrak A}_0),
l=1,2$. So we can apply Theorem~\ref{t:equivalence1}, that
immediately completes the proof of the second part of
Theorem~\ref{main1}. Now Corollary~\ref{main2} follows immediately
from the second part of Theorem~\ref{main1} and
Lemma~\ref{l:vanishing}.

In the case under consideration, we can give a more explicit
description of an operator $U$ that provides Murray-von Neumann
equivalence of the projections ${\mathcal
V}_1E_1(\lambda_1){\mathcal V}_1^{-1}$ and ${\mathcal
V}_2E_2(\lambda_2){\mathcal V}_2^{-1}$. This is based on the
following

\begin{lemma}\label{l1}
For any $a\in {\mathfrak A}_0=\B(\Gamma, \sigma) $, the bounded
operator ${\mathcal V}_2i_2Jp_1{\mathcal V}_1^{-1}\pi(a)$ in
$\ell^2(\Gamma)\otimes \ell^2(\N)$ belongs to $\pi({\mathfrak
A}_0)$.
\end{lemma}

\begin{proof}
We have ${\mathcal V}_2i_2Jp_1{\mathcal V}_1^{-1}=\id\otimes
V_2j_2J_0r_1V_1^{-1}$. It follows that, for any $a\in
C^*_r(\Gamma, \bar\sigma)\otimes {\mathcal K}$, the operator
${\mathcal V}_2i_2Jp_1{\mathcal V}_1^{-1}\pi(a)$ belongs to
$\pi(C^*_r(\Gamma, \bar\sigma)\otimes {\mathcal K})$. Moreover,
the bounded operator ${\mathcal V}_2i_2Jp_1{\mathcal V}_1^{-1}$ in
$\ell^2(\Gamma)\otimes \ell^2(\mathbb N)$ commutes with $D\otimes
\id $ that implies $$ {\mathcal V}_2i_2Jp_1{\mathcal V}_1^{-1}\in
\bigcap_{k\in \mathbb N} \Dom \widetilde\partial^k. $$ Since the
space $\{P\in \bigcap_{k\in \mathbb N} \Dom \widetilde\partial^k :
\widetilde\partial^k(P)\circ (\id\otimes \Delta)\ \text{is
bounded}\ \forall k\in\N\}$ is a left ideal in $\bigcap_{k\in
\mathbb N} \Dom \widetilde\partial^k$ (cf. (\ref{e:ideal})), we
get ${\mathcal V}_2i_2Jp_1{\mathcal V}_1^{-1}\pi(a)\in
\B(\Gamma,\sigma)$ for any $a\in \B(\Gamma,\sigma)$ as desired.
\end{proof}

By Lemma~\ref{l1}, the operator $T={\mathcal
V}_2E_2(\lambda_2)i_2Jp_1E_1(\lambda_1){\mathcal V}_1^{-1}$ in
$\chh$ belongs to $\pi(\mathfrak{A}_0)$, and the operator $U$ that
provides Murray-von Neumann equivalence of the projections
${\mathcal V}_1E_1(\lambda_1){\mathcal V}_1^{-1}$ and ${\mathcal
V}_2E_2(\lambda_2){\mathcal V}_2^{-1}$ can be taken from the polar
decomposition of this operator as in the proof of
Theorem~\ref{t:equivalence1}.

\section{Quantum Hall effect}\label{Hall}

Since the results that we have obtained were essentially known for
Euclidean space and applied to the Euclidean space model for the
integer quantum Hall effect, we will focus on the hyperbolic space
model for the fractional quantum Hall effect, \cite{CHMM},
\cite{MM}.

\subsection{The Hamiltonian}

We begin by reviewing the construction of the Hamiltonian. First
we take as our principal model of hyperbolic space, the hyperbolic
plane. This is the upper half-plane $\hyp$ in $\comp$ equipped
with its usual Poincar\'e metric $(dx^2+dy^2)/y^2$, and symplectic
area form $\omega_\hyp  = dx\wedge dy/y^2$. The group $\slr$ acts
transitively on $\hyp$ by M\"obius transformations
$$
x+iy = \zeta
\mapsto g.\zeta = \frac{a\zeta+b}{c\zeta+d},\quad\mbox{for }
g=\left(
\begin{array}{cc}
a & b\\
c & d
\end{array}\right).$$
Any Riemann surface of genus $g$ greater than 1 can be realised as
the quotient of
$\hyp$ by the action of its fundamental group realised as a
discrete subgroup $\Gamma$ of $\slr$.

Choose a 1-form $\bf A$ called a vector potential whose curvature
is the uniform magnetic field ${\bf B}= d{\bf A} =
\theta\omega_\hyp$, whose flux is $\theta$. As in geometric
quantisation we may regard $\bf A$ as defining a Hermitian
connection $\nabla = d+i{\bf A}$ on the trivial line bundle $\cl$
over $\hyp$, whose curvature is $i\theta\omega_\hyp$. Using the
Riemannian metric the Hamiltonian of an electron in this field is
given in suitable units by $$H = H_{{\bf A},V} =\nabla^*\nabla +
\mu^{-2}V = (d+i{\bf A})^*(d+i{\bf A}) + \mu^{-2}V,$$ where $V$ is
an electric potential associated to a real material and $\mu$ is
the coupling constant. $V$ is also assumed to be invariant under
$\Gamma$ respecting a crystalline type structure. It can be
checked that $H$ commutes with the projective $(\Gamma,
\sigma)$-action on $L^2(\mathbb H)$ as defined in the earlier
sections.

\subsection{Algebra of observables}
Let $\mathcal{F}$ be a connected fundamental domain for the action
of $\Gamma$ on $\hyp$. Take any second order self-adjoint
$(\Gamma,\sigma)$-invariant elliptic differential operator $P$ in
$L^2(\hyp)$, for instance, $P=H_{{\bf A},V}(1)$. Define the
unitary isomorphism $V_2:L^2(\F)\to \ell^2(\N)$, using the
orthonormal basis $\{\varphi_j:j\in\N\}$ of eigenvectors for the
operator $P$ in $L^2(\F)\cong L^2(\hyp/\Gamma)$ with the
$(\Gamma,\sigma)$-periodic boundary conditions:
$V_2\varphi_j=\delta_j, j\in\N$. Introduce a unitary operator
${\mathcal V}_2 : L^2(\hyp)\to \ell^2(\Gamma) \otimes\ell^2(\N)$
as ${\mathcal V}_2=(\id\otimes V_2)\circ {\bf U}$, where ${\bf U}
: L^2(\hyp)\cong \ell^2(\Gamma)\otimes L^2(\mathcal{F})$ is the
$(\Gamma, \sigma)$-equivariant isometry \eqref{isom}. This
operator induces an isomorphism of the algebra ${\mathcal
U}_{L^2(\hyp)}(\Gamma,\bar\sigma)$, consisting of operators on
$L^2(\hyp)$ that commute with the projective $(\Gamma,
\sigma)$-action, with the von Neumann algebra
$\A^L_{\ell^2(\N)}(\Gamma,\bar\sigma)\cong
\A^L(\Gamma,\bar\sigma)\otimes {\mathcal B}(\ell^2(\N))$.

Define the {\em algebra of observables} to be ${\mathcal
B}(\Gamma, \sigma)$ introduced at the beginning of section
\ref{smooth}, which is considered as a $\ast$-subalgebra of
${\mathcal U}_{L^2(\hyp)}(\Gamma,\bar\sigma)$. Recall that we have
established in section \ref{smooth} that $e^{-tH} \in {\mathcal
B}(\Gamma, \sigma)$. The observables of the model include those
spectral projections of the Hamiltonian $H$ corresponding to gaps
in the spectrum. The fact that such a projection belongs to
${\mathcal B}(\Gamma, \sigma)$ was established in section
\ref{smooth} by using the Riesz representation for the projection
and the fact that ${\mathcal B}(\Gamma, \sigma)$ is closed under
the holomorphic functional calculus. This justifies the choice of
${\mathcal B}(\Gamma, \sigma)$  as the algebra of observables.

\subsection{Canonical derivations on the algebra of observables }

Let $\Sigma_g = \hyp/\Gamma$ be the Riemann surface determined by
quotienting by $\Gamma$.
We follow the usual conventions (see for example \cite{GH}) in
fixing representative
homology generators corresponding to cycles $A_j,B_j, j=1,2,\ldots,g$
with each pair $A_j, B_j$ intersecting in a common base point and
all other intersection numbers being zero. Let
$a_j, j=1,\ldots,g$ be harmonic 1-forms dual to $A_j, j=1,\ldots,g$
and $b_{j}, j=1,2,\ldots,g$ be harmonic 1-forms dual to
$B_j,j=1,2,\ldots,g$. Let $\widetilde a_j, \widetilde b_j$
denote the lifts of $a_j, b_j$ to $\mathbb H$ respectively.

Define the functions on $\mathbb H$ given by,
$$
\Omega_j(z) = \displaystyle i\int_u^z \widetilde a_j, \qquad
\Omega_{j+g}(z) =
\displaystyle i\int_u^z \widetilde b_j,
$$
where $u\in\mathbb H$ is a fixed point. Since $\widetilde a_j, \widetilde b_j$
are bounded 1-forms on $\mathbb H$, one sees that
there are positive constants $C_j$ such that
\begin{equation}\label{est}
|\Omega_j(z) |\le C_j d(u, z) \qquad \text{for all} \; z\in \mathbb H.
\end{equation}
For any $j=1,2,\ldots, 2g$, denote by $\Omega_j$ the operator in
$L^2(\hyp)$ of multiplication by the function $\Omega_j$. Define
$$\delta_j(T) = [\Omega_j,T],$$ where $T\in \C(\Gamma,\bar\sigma)
\otimes {\mathcal R}\subset {\mathcal B}(\Gamma,\bar\sigma)$ is
considered as a bounded operator in $L^2(\hyp)$.

\begin{lemma}\label{deltaj}
For any $T\in \C(\Gamma,\bar\sigma)\otimes {\mathcal R}$, the
operator $\delta_j(T)$ is in $\C(\Gamma,\bar\sigma)\otimes
{\mathcal R}$.
\end{lemma}

\begin{proof} Using the arguments given in the proof of Lemma~\ref{F},
one can easily see that, under the isomorphism $T\in
\A^L_{\ell^2(\N)}(\Gamma,\bar\sigma) \mapsto {\mathcal V}^*_2 T
{\mathcal V}_2\in {\mathcal U}_{L^2(\hyp)}(\Gamma,\bar\sigma),$
the $\ast$-subalgebra $\C(\Gamma,\bar\sigma)\otimes {\mathcal
R}\subset \A^L_{\ell^2(\N)}(\Gamma,\bar\sigma)$ corresponds to the
subalgebra in ${\mathcal U}_{L^2(\hyp)}(\Gamma,\bar\sigma)$, which
consists of all $(\Gamma,\sigma)$-invariant bounded operators $Q$
in $L^2(\hyp)$ whose Schwarz kernels $k_Q$ are smooth and properly
supported (i.e. $k_Q(x,y)=0$ when $d(x,y)>C$ with some constant
$C>0$). Recall that $(\Gamma,\sigma)$-invariance of $Q$ is
equivalent to the relation (\ref{kernel}) for its Schwarz kernel
$k_Q$.

Take any $T\in \C(\Gamma,\bar\sigma)\otimes {\mathcal R}$
considered an a bounded operator in $L^2(\hyp)$. Let $k_T$ be its
Schwarz kernel. Then the Schwarz kernel $k_{[\Omega_j,T]}$ of the
operator $[\Omega_j,T]$ is given by
\[
k_{[\Omega_j,T]}(x,y)=(\Omega_j(x)-\Omega_j(y))k_T(x,y), \quad
x,y\in \M.
\]
Clearly, $k_{[\Omega_j,T]}$ is smooth and properly supported.

It is easy to see that, for $\gamma\in \Gamma$, the difference $
\Omega_j(\gamma.z) - \Omega_j(z)$ is constant independent of
$z\in\hyp$. Using this fact, one can check that $k_{[\Omega_j,T]}$
satisfies (\ref{kernel}), that completes the proof.
\end{proof}

Therefore $\delta_j, j=1,2,\ldots, 2g$, is a (densely defined)
derivation on the algebra of observables $\ck(\Gamma, \sigma)$.

\subsection{Hall conductance cyclic cocycle}

In this subsection we recall the Kubo formula for the
          Hall conductivity.
The reasoning is that
          the Hall conductivity
is measured  by determining the equilibrium
          ratio of the current in the direction of the applied electric
field to the
          Hall voltage, which is the potential difference in the orthogonal
direction.
          To calculate this mathematically we instead determine
          the component of the induced current that is orthogonal to the applied
          potential. The conductivity can then be obtained by dividing this
          quantity by the magnitude of the applied field.
Interpreting the generators of the fundamental group as geodesics
on hyperbolic space gives a family of preferred directions
emanating from the base point. One of the basic results in
\cite{CHMM} is the following: the expectation of the current $J_k$
is given by $$\Tr_\Gamma(P\delta_kH) =
i\Tr_\Gamma(P[\partial_tP,\delta_kP]) =
-iE_j\Tr_\Gamma(P[\delta_jP,\delta_kP]),$$  where $E_j$ is the
electric field in the $j$ direction. Therefore one sees that,

{\em The conductivity for currents in the $k$ direction induced by
electric fields in the $j$ direction is given by
$-i\Tr_\Gamma(P[\delta_jP,\delta_kP])$. }

The following is Lemma 12 in \cite{CHMM}.

\begin{lemma}
For any $j,k=1,\ldots,2g$, the formula $$c_{j,k}(T_0, T_1, T_2)=
\Tr_\Gamma(T_0[\delta_jT_1,\delta_kT_2]) =
\Tr_\Gamma(T_0[\Omega_j, T_1][\Omega_k,T_2]), \quad T_0, T_1, T_2
\in  \C(\Gamma,\bar\sigma)\otimes {\mathcal R}, $$ defines a
cyclic $2$-cocycle on $\C(\Gamma,\bar\sigma)\otimes {\mathcal R}$.
\end{lemma}

\begin{defn*} The Kubo formula for the Hall conductance cyclic
$2$-cocycle $\tr_K$ on a dense subalgebra
$\C(\Gamma,\bar\sigma)\otimes {\mathcal R}$ of the algebra of
observables  $\ck(\Gamma, \sigma)$ is defined as $$\tr_K =
\sum_{j=1}^g c_{j,j+g}.$$
\end{defn*}

There is a symplectic map from $\hyp$ to $\real^{2g}$ given by
$\Xi:z\mapsto (\Omega_1(z),\ldots,\Omega_{2g}(z))$. It is the lift
to $\hyp$ of the Abel-Jacobi map on $\mathbb H/\Gamma$. Define a
group $2$-cocycle $\Psi: \Gamma\times\Gamma\to \mathbb R$ as
follows. Consider the straight-edged triangle $\Delta(u, \gamma_1,
\gamma_2)$ in $\real^{2g}$ obtained  by joining the 3 points
$\Xi(u)=0$, $\Xi(\gamma_1^{-1}.u)$ and $\Xi(\gamma_2.u)$. Then
$\Psi(\gamma_1, \gamma_2)$ is defined to be the symplectic area of
$\Delta(u, \gamma_1, \gamma_2)$, which is equal to
$$\sum_{j=1}^g(\Omega_j(\gamma_1^{-1}.u)\Omega_{j+g}(\gamma_2.u)-
\Omega_{j+g}(\gamma_1^{-1}.u)\Omega_{j}(\gamma_2.u)).$$

We have seen earlier that $|\Omega_{j}(\gamma.u)|\le C_j\; d(u,
\gamma.u) \le C_j' \ell(\gamma)$. Using the Cauchy-Schwartz
inequality and the fact that $\ell(\gamma) =\ell(\gamma^{-1})$, we
see that there is a positive constant $C$ such that $$
|\Psi(\gamma_1, \gamma_2)|\le C (1+ \ell(\gamma_1))^2 (1 +
\ell(\gamma_2))^2\qquad \rm{for \; all}\;\;\gamma_1, \gamma_2\in
\Gamma. $$ That is, $\Psi$ is a polynomially bounded group
2-cocycle on $\Gamma$. Recall that the group $2$-cocycle $\Psi$
defines  the cyclic $2$-cocycle $\tau_\Psi$ on
$\C(\Gamma,\bar\sigma)$. By Lemma \ref{RD}, we see that
$\tau_\Psi\#\Tr$ extends by continuity to $\B(\Gamma, \sigma)$.

One of the main results of \cite{CHMM} is the following.

\begin{thm}
The Hall conductance cyclic $2$-cocycle $\tr_K$ agrees with
$\tau_\Psi\#\Tr$ on $\C(\Gamma,\bar\sigma)\otimes {\mathcal R}$.
\end{thm}

Hence by Lemma \ref{RD}, the Hall conductance cyclic
$2$-cocycle $\tr_K$ also extends by continuity to $\B(\Gamma, \sigma)$.
Therefore we are in a position to apply the Corollary \ref{main2}
to deduce Corollary \ref{main4}.


\begin{thebibliography}{Comtet+Hs}


\bibitem[Bel]{Bel}{J. Bellissard, A. van Elst, H. Schulz-Baldes}, {The
non-commutative geometry of the quantum Hall effect}, {\em J. Math.
Phys.} {\bf 35} (1994),
{5373-5451}.

\bibitem[Bost]{Bost} J. Bost, Principe d'Oka, $K$-th\'eorie et syst\`emes
dynamiques
non commutatifs. {\em Invent. Math.} {\bf 101} (1990) 261-333.

\bibitem[BrSu]{BrSu} J. Br\"uning, T. Sunada, On the spectrum of
gauge-periodic elliptic operators. M\'ethodes semi-classiques,
Vol. 2 (Nantes, 1991). {\em Ast\'erisque} {\bf 210} (1992), 65-74.

\bibitem[Bu]{Bu} U. Bunke,
On the gluing problem for the $\eta$-invariant.
{\em J. Differential Geom.}, \textbf{41} (1995), 397--448.

\bibitem[BFKM]{BFKM}
D. Burghelea, L. Friedlander,T. Kappeler, P. McDonald,
Analytic and Reidemeister torsion for representations in finite type 
Hilbert modules.
{\em Geom. Funct. Anal.} \textbf{6} (1996), 751--859.

\bibitem[CHMM]{CHMM} A. Carey, K. Hannabus, V. Mathai and P. McCann,
Quantum Hall Effect on
the hyperbolic plane, {\em Commun. Math. Physics}, {\bf 190} no. 3
(1998) 629-673.

\bibitem[Co81]{Co81} A. Connes, An analogue of the Thom
isomorphism for crossed products of a $C^*$ algebra by an action
of $\mathbb R$, {\em Adv. in Math.} {\bf 39} (1981), 31-55.

\bibitem[Co]{Co} A. Connes, Noncommutative differential geometry,
{\em Publ.Math.I.H.E.S.} {\bf 62} (1986), 257-360.

\bibitem[Co2]{Co2} A. Connes, {\em Noncommutative geometry}. Acad.
Press, Inc., San Diego, CA,
(1994).

\bibitem[CM90]{CM90} A. Connes, H. Moscovici, Cyclic cohomology,
the Novikov conjecture and hyperbolic groups, {\em Topology} {\bf 29}
(1990), 345-388.

\bibitem[G]{Greiner} P. Greiner, An asymptotic expansion for the heat equation,
{\em Arch. Ration. Mech. and Anal.}, {\bf 41} (1971), 163-218.

\bibitem[GH]{GH} P. Griffiths and J.  Harris, {\em Principles of algebraic
geometry,} Pure and Applied Mathematics. Wiley-Interscience [John
Wiley \& Sons], New York, 1978.

\bibitem[Gr]{Gr} M. Gromov, Volume and bounded cohomology, {\em Publ.
Math. I.H.E.S.}
{\bf 56} (1982), 5-99.

\bibitem[Gr2]{Gr2} M. Gromov, Hyperbolic groups. Essays in group
theory, 75--263, Math. Sci. Res. Inst.
Publ., 8, Springer, New York, 1987.

\bibitem[Gui]{Gui} A. Guichardet, {\em Cohomologie des groupes topologiques
et des alg\`ebres de Lie.} Cedic/ Fernand Nathan, Paris, 1980.

\bibitem[HS84]{HS84}
B. Helffer, J. Sj\"ostrand,  Multiple wells in the semiclassical limit. I.
{\em Comm. Partial Differential Equations}  \textbf{9}  (1984),  337--408.

\bibitem[HS87]{HS87} B. Helffer, J. Sj\"ostrand,
Effet tunnel pour l'\'equation de Schr\"odinger avec champ magn\'etique.
{\em Ann. Scuola Norm. Sup. Pisa Cl.
Sci.} (4)  \textbf{14}  (1987),  625--657 (1988).

\bibitem[HS88]{HS88} B. Helffer, J. Sj\"ostrand,
Analyse semi-classique pour l'\'equation de Harper (avec application 
\`a  l'\'equation
de Schr\"odinger avec champ magn\'etique).
M\'em. Soc. Math. France (N.S.)  No. 34 (1988).


\bibitem[HS90]{HS90} B. Helffer, J. Sj\"ostrand,
Analyse semi-classique pour l'\'equation de Harper. II. Comportement 
semi-classique
pr\`es d'un rationnel.
M\'em. Soc. Math. France (N.S.)  No. 40 (1990).

\bibitem[HS89]{HS89} B. Helffer, J. Sj\"ostrand,
Semiclassical analysis for Harper's equation. III. Cantor structure 
of the spectrum.
M\'em. Soc. Math. France (N.S.)  No. 39 (1989), 1--124.


\bibitem[HSLNP]{HSLNP} B. Helffer, J. Sj\"ostrand,
\'Equation de Schr\"odinger avec champ magn\'etique et \'equation de Harper.
Schr\"odinger operators (S\o nderborg, 1988),  118--197,
Lecture Notes in Phys., \textbf{345}, Springer, Berlin, 1989.


\bibitem[HSHaas]{HSHaas} B. Helffer, J. Sj\"ostrand,
On diamagnetism and de Haas-van Alphen effect.
Ann. Inst. H. Poincar\'e Phys. Th\'eor.  \textbf{52}  (1990),
303--375.


\bibitem[Ji1]{Ji92} R. Ji, Smooth dense subalgebras of reduced
group $C^*$-algebras, Schwartz cohomology of groups, and cyclic
cohomology, {\em J. Funct. Anal.} {\bf 107} (1992), 1--33.


\bibitem[Ji2]{Ji} R. Ji,  A module structure on cyclic cohomology of
group graded algebras, $K$-Theory {\bf 7} (1993),  369--399.


\bibitem[JiS]{JiS} R. Ji and L. Schweitzer, Spectral invariance of
smooth crossed products,
and rapid decay locally compact groups, {\em $K$-Theory} 10 (1996),
no. 3, 283--305.


\bibitem[Ko]{Ko} Yu. Kordyukov, $L^p$-theory of elliptic differential
operators on manifolds of bounded geometry, {\em Acta Appl.
Math.}, {\bf 23} (1991), {223--260}.


\bibitem[Ko04]{Ko04} Yu. Kordyukov,
Spectral gaps for periodic Schr\"odinger operators with strong magnetic fields.
Preprint math.SP/0311200. (To appear in Commun. Math. Phys.)


\bibitem[MM]{MM} M. Marcolli, V. Mathai, Twisted index theory for
good orbifolds, II:
fractional quantum numbers, {\em Commun. Math. Phys.} {\bf 217}
{(2001)}, {55-87}.


\bibitem[Ma]{Ma} V.~Mathai, On positivity of the Kadison constant and
noncommutative Bloch theory, {\em Tohoku Math. Publ.} {\bf 20}(2001),
{107-124}.


\bibitem[MS]{MS} V. Mathai and M. Shubin, Semiclassical asymptotics
and gaps in the spectra of
magnetic Schr\"odinger operators, {\em Geometriae Dedicata} {\bf 91}(2002),
{155-173}.


\bibitem[Nak+Bel]{Nak+Bel}{S. Nakamura, J. Bellissard}, {Low energy
bands do not contribute
to the quantum Hall effect}, {\em Commun. Math. Phys.} {\bf 131}
{(1990)}, {283-305}.

\bibitem[NovSh]{NovSh}
S. P. Novikov and M.A. Shubin,
Morse inequalities and von Neumann ${\rm II}\sb 1$-factors.
{\em Dokl. Akad. Nauk SSSR} {\bf 289} (1986), no. 2, 289--292.

\bibitem[Schw]{Schw} L. Schweitzer, A short proof that $M_n(A)$ is
local if $A$ is local and Fr\'echet, {\em Internat. J. Math.},
{\bf 3} (1992), 581 - 589.

\bibitem[Sh]{Sh} M. Shubin, {Semiclassical asymptotics on covering
manifolds and
Morse Inequalities}, {\em Geom. Anal. and Func. Anal.}, {\bf 6}, no.
2 (1996), 370-409.

\bibitem[Sh2]{Sh2} M.~Shubin, {Discrete magnetic Laplacian}, {\em
Commun. Math. Phys.}
{\bf 164} (1994), no.2, 259--275.

\bibitem[Ta]{Ta} M. Takesaki,
Theory of operator algebras. I.
Springer-Verlag, New York-Heidelberg, 1979.

\end{thebibliography}
\end{document}